\documentclass[12pt]{amsart}

\usepackage{amsmath,amsfonts,bm}

\def\eqref#1{equation~\ref{#1}}

\def\1{\bm{1}}

\DeclareMathAlphabet{\mathsfit}{\encodingdefault}{\sfdefault}{m}{sl}
\SetMathAlphabet{\mathsfit}{bold}{\encodingdefault}{\sfdefault}{bx}{n}

\newcommand{\R}{\mathbb{R}}

\usepackage{amsmath}
\usepackage{amssymb}
\usepackage{mathtools}

\usepackage{algorithm}
\usepackage[noend]{algpseudocode}

\usepackage{bm}
\usepackage{enumitem}
\usepackage{graphicx}
\usepackage{subcaption}
\usepackage{hyperref}

\usepackage{parskip}
\usepackage{xcolor}

\usepackage{multirow}
\usepackage{booktabs}
\usepackage{tabularray}
\allowdisplaybreaks

\newtheorem{theorem}{Theorem}

\theoremstyle{definition}

\theoremstyle{remark}
\newtheorem{remark}[theorem]{Remark}
\numberwithin{equation}{section}

\newcommand{\Rb}{\mathbb{R}}

\numberwithin{equation}{section}
\renewcommand{\theequation}{\arabic{section}.\arabic{equation}}

\begin{document}
\title[Data-driven optimal control with neural networks] 
{Data-driven optimal control with neural network modeling of gradient flows}
\author{Xuping Tian}
\address{Iowa State University, Department of Mathematics, Ames, IA 50011} 
\email{xupingt@iastate.edu}
\author{Baskar Ganapathysubramanian}
\address{Iowa State University, Department of Mechanical Engineering, Ames, IA 50011} 
\email{baskarg@iastate.edu}
\author{Hailiang Liu}
\address{Iowa State University, Department of Mathematics, Ames, IA 50011} \email{hliu@iastate.edu}
\thanks{This work was supported by the Translational AI Center (TrAC seed grant 2022-2023) at Iowa State University}

\begin{abstract}  
Extracting physical laws from observation data is a central challenge in many diverse areas of science and engineering. We propose Optimal Control Neural Networks (OCN) to learn the laws of vector fields in dynamical systems, with no assumption on their analytical form, given data consisting of sampled trajectories. The OCN framework consists of a neural network representation and an optimal control formulation. We provide error bounds for both the solution and the vector field. The bounds are shown to depend on both the training error and the time step between the observation data. We also demonstrate the effectiveness of OCN, as well as its generalization ability, by testing on several canonical systems, including the chaotic Lorenz system.  
\end{abstract}

\maketitle


\section{Introduction} 

A central challenge in many diverse areas of science 
and engineering is to discover physical laws. This work concerns learning dynamical systems arising from real-world applications but where a complete mathematical description of the dynamics is unavailable. In such scenarios, we rely on extracting insight from data. 
Our work sits at the intersection of machine learning and dynamical systems, where the equations describing the dynamics are implicitly reconstructed from observed trajectory data using neural networks.    

{\bf Data-driven discovery of dynamical systems.}
There is a long and fruitful history of modeling dynamics from data. Earlier efforts for system discovery include a large set of methods (See Section \ref{rwks} below). One fruitful family of approaches includes using symbolic regression \cite{BL07, SL09} for finding nonlinear equations. This strategy balances the complexity of the model with predictive power. These approaches are often expensive and require careful selection of candidate models or basis expressions.  
More recently, sparsity has been used to determine the governing dynamical system \cite{BP+16, CL+19, RB+17, S17, ZA+18}, where certain sparsity-promoting strategies are deployed to obtain parsimonious models. 
The challenge with this strategy lies in choosing a suitable sparsifying function basis.  
There also have been studies on system identification using Gaussian processes \cite{KGBM05,RK18} and statistical learning \cite{LZT19}. Instead of discovering the exact function(al) expressions, one also seeks to reconstruct accurate numerical approximations to the dynamical systems; see e.g. \cite{RP+18, R18, QW+19, LJ+19, LL19, DGYZ22, SZS20, LM+22} for works using the neural network representation. Our work in this paper falls into the latter category.


{\bf Deep neural networks (DNN).} DNNs have seen tremendous success in many disciplines, particularly supervised learning. Their structure with numerous consecutive layers of artificial neurons allows DNNs to express complex input-output relationships. Efforts have been devoted to the use of DNNs for various aspects of scientific computing, including solving and learning systems involving ODEs and PDEs. Recently, the interpretation of residual networks by He et al. \cite{He16} as approximate ODE solvers in \cite{E17} spurred research on the use of ODEs in deep learning \cite{CR+18, LZ+18, HR17}. Neural ODEs \cite{CR+18} as neural network models generalize standard layer-to-layer propagation to continuous depth models. Along this line of research, work \cite{LM20} develops a PDE model to represent a continuum limit of neural networks in both depth and width. 

{\bf Optimal control neural networks.} Recently, there has been a growing interest in understanding deep learning methods through the lens of dynamical systems and optimal control \cite{LC+18,LH18,ZZ+19,BC+19}. An appealing feature of this approach is that the compositional structure is explicitly taken into account in the time evolution of the dynamical systems, from which novel algorithms and network structures can be designed using optimal control techniques. 
In particular, mathematical concepts from optimal control theory are naturally translatable to dynamic neural networks, and provide interesting possibilities, including computing loss gradient by the adjoint method and natural incorporation of regularization and/or prior knowledge into the loss function. This work directly takes advantage of these concepts. 

In this paper, we build upon recent efforts that discover dynamical systems using deep neural networks (DNNs) \cite{RP+18, QW+19} and the optimal control approach for learning system parameters \cite{LT21}. We seek to gain new insight into the dynamics discovery problem using ``optimal control networks'' (OCN for short). Taking gradient flows $\dot x =-\nabla f(x)$ as a model class, we establish mathematically sound, dynamically accurate, computationally efficient techniques for discovering $f$ from trajectory data. Note that the values of $f$ are not observed, in contrast to the standard supervised learning problems. We exploit the representation power of deep neural networks to approximate $f$, unlike related recent efforts that require feature libraries \cite{BP+16, CL+19, RB+17, S17, ZA+18}. The key steps involved in OCN include:  
{ 
\begin{enumerate}
\item We exploit a neural network $G(\cdot,\theta)$ as a global representation of the unknown governing function $f$, where $\theta$ represents the neural network parameters to be learned.
\item We then formulate the learning problem as an optimal control problem of form  
\begin{align*}
  \min_{\theta\in \mathcal{A}} \quad &J(\theta) = \sum_{i=1}^{n}L_i(y(t_i)),\\ 
  \text{s.t.}\quad &\dot y(t)=-\partial_y G(y(t),\theta) \quad t\in(t_0,T], \quad y(t_0)=x_0,
\end{align*}    
where $\mathcal{A}\subset \Rb^N$ is the control set, $t_n=T$ and
\begin{align*}
  L_i(y) &:= \|y-x_i\|^2,\quad 1\le i\le n.
\end{align*}
Here $x_i$ is the observed data at time $t_i$, $L_i$ is a local loss that measures the error between the solution to ODE in the constraint and the observed data at $t_i$.
\item We apply a gradient-based method to update the network parameters $\theta$, where the loss gradient is evaluated by
\begin{equation*}
\nabla_\theta J= -\sum_{i=0}^{n-1}\int_{t_{i}}^{t_{i+1}}   \big(\partial_\theta \partial_y G(y(t),\theta)\big)^\top p(t)dt.   
\end{equation*}
Here, both the state variable $y$ and the co-state variable $p$ are obtained by solving the coupled system: 
\begin{align*}
\dot y(t) & =-\partial_y G(y(t), \theta), \quad y(t_0)=x_0, \\
\dot p(t)&=\big(\partial^2_y G(y(t),\theta)\big)^\top p(t), \quad t_{i-1}\le t<t_i,\quad i=n,...,1,\\
  p(T) &=\partial_y L_n(y(T)),\;
  p(t_i^-) = p(t_i^+) + \partial_y L_i(y(t_i)), \quad i=n-1,...,1. 
\end{align*}
\item In order to achieve high-order accuracy of the gradient evaluation in (3), we apply a partitioned Runge-Kutta method to solve the coupled $(y, p)$ system. The Runge-Kutta solver is shown to be symplectic in the sense that it conserves the bilinear quantity $\Big(\frac{\partial y(t)}{\partial y(t_0)}\Big)^\top p(t)$ for $t\in (t_{i-1}, t_i]$ where $i=1, \cdots, n$. This is crucial since such a bilinear quantity is an invariant of the continuous system.
\end{enumerate}
}
The methodology and key formulations apply directly to more general dynamical systems $\dot x = F(x)$ and can be generalized to parameterized, time-varying, or externally forced systems. 

This paper makes the following specific contributions:
\begin{itemize}
    \item We propose and analyze a novel framework for discovering dynamical systems from the observation data, incorporating neural network approximations into an optimal control formulation. 
    \item We establish error bounds for both the solution and the vector field, which show that the global error depends only on the training error and the time step between the observation data.
    \item {We incorporate a partitioned symplectic Runge-Kutta method into the training process of the OCN neural network,  which is a symplectic solver and guarantees a high-order accuracy of the loss gradient estimation. }  
    \item We demonstrate the effectiveness and generalization ability of OCN on several canonical systems. In particular, we provide a thorough exploration on the chaotic Lorenz system, which suggests that OCN exhibits superior performance (to symbolic approaches like SINDy \cite{BP+16}) when the derivative data $\dot x$ is unavailable or the data $x$ has relatively large time steps.
\end{itemize}

\subsection{Further related works}\label{rwks}
There are techniques that address various aspects of 
the dynamical system discovery problem, including methods to discover governing equations from time series data \cite{CM87}, equation-free modeling \cite{KG+03}, empirical dynamic modeling \cite{SM12, YB+15},  modeling emergent behavior \cite{R14}, nonlinear Laplacian spectral analysis \cite{GM12}, artificial neural networks \cite{GR+98}, Koopman analysis \cite{WK+15, BB+17, AKM18}, learning the effective dynamics \cite{VAUK20, VZPK21} and automated inference of dynamics \cite{DN15, SV+11}. Instead of reconstructing the dynamical systems, there are also works that focus on learning the parameters in some dynamical systems  \cite{AF+18,LE+22,LT21}.

{\bf Training of neural ODEs.} This work is also complementary to efforts that incorporate ODE solvers into training neural networks, including numerical methods for training neural ODEs. 
Using the adjoint method to train neural networks was first introduced in \cite{CR+18}. To overcome the numerical errors associated with this approach, several techniques have been proposed, for instance, the checkpoint method \cite{GK19, ZD+20},  the asynchronous leapfrog method \cite{ZD+21}, the symplectic adjoint method \cite{MM21}, interpolation method \cite{DK+20}, and the proximal implicit solvers \cite{BX+22}. 


  


{\bf Structure-preserving learning.} For many application problems, it is desirable to adopt structured machine learning approaches, where one imparts from the outset some physically motivated constraints to the dynamical system to be learned. The gradient flow dynamics learned by our approach have a precise physical structure, which not only ensures the stability of the learned models automatically, but also gives physically interpretable quantities.  Such advantages have been observed by researchers when learning different structural systems, such as stable dynamic systems \cite{KM19,GHRW20}, Hamiltonian systems \cite{GDY19, ZDC20a, JZKT20, BDMK20}, and more general systems based on a generalized Onsager principle \cite{YTEL21}.  


Our work aligns with \cite{RP+18, QW+19} but with a different strategy. 
Work in \cite{QW+19} first discretizes the dynamical system based on a local integral form, then uses a neural network to approximate the local flow map between two neighboring data points.  Such strategy may be seen more as learning of an ODE solver specified through the loss function.  In contrast, we incorporate a global network representation into the optimal control formulation.   
Such global approximation using neural network representation is also considered in \cite{RP+18}, however, the parameter learning method therein is built for a discretized dynamical system in the form of multi-step time-stepping schemes. Importantly, we are able to obtain error bounds that allow users to judiciously reason about the accuracy and convergence of our method. 

The rest of the paper is arranged as follows: problem setup and our method are introduced in Section 2 with detailed mathematical formulations. Section 3 presents a theoretical analysis of the errors. Computational details of our method are presented in Section 4. Section 5 includes several numerical experiments. Finally, some concluding remarks and discussions are given in Section 6. Implementation details and technical proofs are given in the appendix. 


\section{Method}
Here we provide an overview of our method. We first present the problem setup based on a set of time series data in order to learn the unknown vector field.  Afterwards, we argue why we can use neural networks to realize the needed approximation. Then we explain the learning phase of the neural network, which seeks to solve an optimal control problem. Finally, we explain the training stage, where we are able to produce gradients in parameter space to update network parameters.

\begin{figure}[t]
\centering
\includegraphics[width=0.8\linewidth]{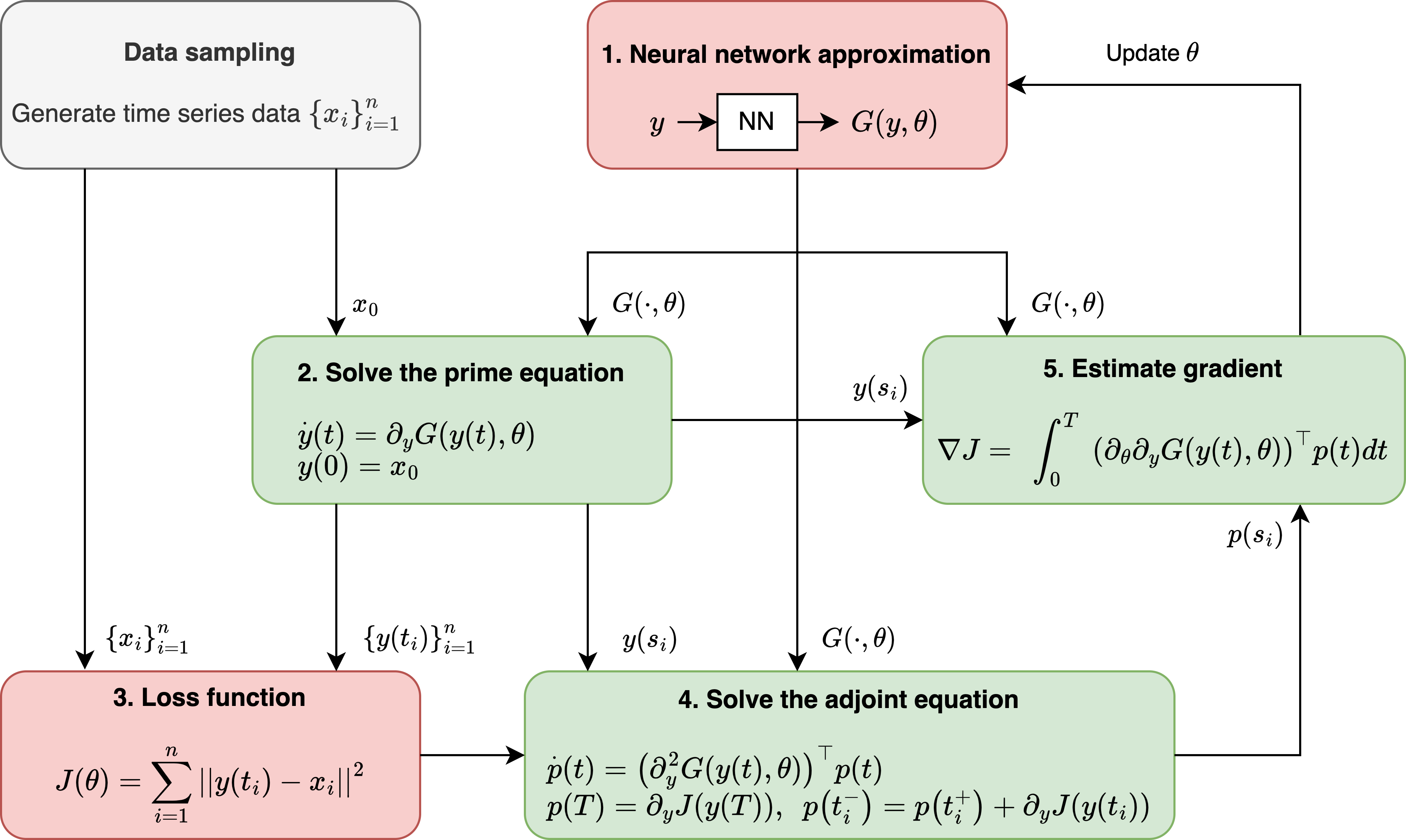}
\caption{A concept diagram showing the flow of OCN. Steps 2, 4, and 5 correspond to the gradient evaluation method presented in Theorem \ref{thm1}.}
\label{fig:ocn}
\end{figure}

\subsection{Problem setup}\label{psu} 
Many application problems are modeled by gradient flows \cite{AGS05}. We consider gradient flow systems of the form
\begin{equation}\label{gf}
\dot x(t)=-\nabla f(x(t)),\quad x(0)=x_0 
\end{equation}
on $[0, T]$, where $x\in\Rb^d$ is the state variable. 
In this paper, we assume the form of $f:\Rb^d\to\Rb$ is unknown. We aim to create an accurate model for 
learning or recovering $f$ using data sampled from solution trajectories and generating solutions over a specified time interval. 

Numerically, in order to produce  trajectories  of  the  dynamical  system when $f$ is known, one  can use various integrators, such as forward Euler,
\begin{align}\label{fe}
x_{i+1} =x_i- \Delta t \nabla f(x_i),
\end{align}
where the time domain is divided into equal step sizes $\Delta t$ so that $t_{i+1}-t_i=\Delta t$ for $0=t_0 <...<t_n=T$. Other high-order accuracy schemes e.g.  4th order Runge-Kutta can also be used. Here we assume that data is collected as solution states on a uniform lattice of time points $\{t_i\}_{i=0}^n$.

\subsection{Neural network approximation} 
The universal approximation theorem states that any continuous function can be approximated arbitrarily well by a neural network \cite{HS+89, B93}. We therefore choose to represent $f(x)$ using a neural network. 

A fully connected feedforward neural network $G(\cdot,\theta):\Rb^{N_1}\to \Rb^{N_m}$ can be seen as a composition of a sequence of linear functions and nonlinear functions:
\begin{equation}\label{Gnn}
G(\cdot,\theta) := \sigma_{m-1} \circ h_{m-1} \circ \dots \circ \sigma_1 \circ h_1.    
\end{equation}
Here $h_j: \Rb^{N_j}\to\Rb^{N_{j+1}}$ are linear functions:
$
h_j(x) = W_j x+b_j,
$
where $W_j\in\Rb^{N_j\times N_{j+1}}$ are matrices, also called weights, $b_j\in\Rb^{N_{j+1}}$ are biases. $\sigma_j:\Rb\to\Rb$ are nonlinear activation functions applied component-wisely to the $j$-th layer. $\theta\in \Rb^N$ denotes the parameter set containing all the parameters $W_1, b_1,...,W_{m-1}, b_{m-1}$ involved, where $N=\sum_{j=1}^{m-1}(N_j+1)N_{j+1}$. 
Some common choices for the activation functions are hyperbolic tangent functions, sigmoid functions, ReLU, etc. \cite{AD+21}.




\subsection{Loss function} 
Though our goal is to learn the function $f$, with no access to function values $f(x_i)$, the usual supervised learning is not applicable.  
The way we learn the parameter $\theta$ of the neural network $G$ is to solve the parameterized ODE system 
\begin{equation}\label{ode-nn}
\dot y(t)=-\partial_y G(y(t),\theta), \quad y(0)=x_0,
\end{equation}
and compare the solution at $t_i$ with the observed data $x_i$. To this end, we take the loss function  
 \begin{equation}\label{loss}
J(\theta) = \sum_{i=1}^{n}\|y(t_i)-x_i\|^2,    
\end{equation}
where the dependence of $J$ on $\theta$ is through $y(t)$. 

In this work, we focus on (\ref{gf}), which is an autonomous system, i.e. $f$ depends solely on the state variable $x$, but not on time $t$, hence $\theta$ can be a time-independent parameter. This point is important for our choice of numerical solvers for (\ref{ode-nn}).

\subsection{Optimal control formulation} Now our problem is reduced to learning $\theta$ 
 by minimizing the loss function (\ref{loss}) subjected to the ODE system (\ref{ode-nn}). 
 From the perspective of control, we need to find an optimal parameter $\theta^*$ for (\ref{ode-nn}) such that the loss function (\ref{loss}) is minimized. This motivates us to formulate it as an optimal control problem:
\begin{subequations}\label{opcr}
\begin{align}
  \min_{\theta\in \mathcal{A}} \quad &J(\theta) = \sum_{i=1}^{n}L_i(y(t_i)),\\ 
  \text{s.t.}\quad &\dot y(t)=-\partial_y G(y(t),\theta) \quad t\in(0,T], \quad y(0)=x_0,
\end{align}    
\end{subequations}
where $\mathcal{A}\subset \Rb^N$ is a control set, $t_n=T$.
Here $L_i$ is a local loss that measures the error between the solution to (\ref{opcr}b) when $y=y(t_i)$ and the observed data $x_i$ at $t_i$. When $n=1$, this reduces to the usual optimal control with terminal cost. We solve this optimal control problem by iteration with gradient-based methods to update $\theta$. For instance, given $\theta_k$, gradient descent (GD) computes $\theta_{k+1}$ by
\begin{equation}\label{gd}
\theta_{k+1} = \theta_k - \eta_k\nabla J(\theta_k), 
\end{equation}
where $\eta_k$ is the step size. One of the main tasks here is to compute the gradient  $\nabla J(\theta)$. This can be obtained via backpropagation through ODE solvers, which gives a discrete approximation to the dynamical system. Another approach to computing the gradient is to use the adjoint method, which is summarized in Theorem \ref{thm1}. 

\subsection{Compute the gradient}
 The following result allows computation of the gradient $\nabla J(\theta)$.

\begin{theorem}\label{thm1} For problem (\ref{opcr}), if $(y(t), \theta)$, $0\leq t\leq T$ is the state trajectory starting from $x_0$, then there exists a co-state trajectory $p(t)$ satisfying 
\begin{subequations}\label{00}
\begin{align}
\dot y(t) & =-\partial_y G(y(t), \theta), \quad y(0)=x_0,\\
\dot p(t)&=\big(\partial^2_y G(y(t),\theta)\big)^\top p(t), \quad t_{i-1}\le t<t_i,\quad i=n,...,1,\\
  p(T) &=\partial_y L_n(y(T)),\;
  p(t_i^-) = p(t_i^+) + \partial_y L_i(y(t_i)), \quad i=n-1,...,1. 
\end{align}
\end{subequations}
Moreover, the gradient of $J$ can be evaluated by 
\begin{equation}\label{grad-}
\nabla J= -\sum_{i=0}^{n-1}\int_{t_{i}}^{t_{i+1}}   \big(\partial_\theta \partial_y G(y(t),\theta)\big)^\top p(t)dt.   
\end{equation}
\end{theorem} 
This allows us to compute $\nabla J$ at each iteration, say when $\theta=\theta_k$, in three steps:
\begin{enumerate}
\item[Step 1.] Solve the forward problem to obtain state $y_k(t):=y(t; \theta_k)$,
\item[Step 2.] Solve the piece-wise backward problem to obtain co-state $p_k(t)$,
\item[Step 3.] Evaluate the gradient of $J$ by (\ref{grad-}), which gives the needed 
$\nabla J(\theta_k)$.
\end{enumerate}
{We shall discuss the computational procedure for the adjoint method in Section \ref{discre}.}



In practice, some real-world systems are not in the form of gradient flows, and our framework is readily extended to encompass these situations, allowing for the discovery of general ODE systems
\begin{equation}
\dot x(t)=F(x(t)),
\end{equation}
where $F:\Rb^d\to\Rb^d$ is unknown. In such case, Theorem \ref{thm1} needs to be modified by replacing $-\partial_yG(y(t),\theta)$ with $G(y(t),\theta)$, where $G(\cdot,\theta):\Rb^d\to\Rb^d$ is a neural network approximator of $F$. We also conducted some numerical tests on this type of problem; see Section \ref{nodetest}. Finally, we should point out that \textit{any priori} knowledge of the properties of $G$  could be used to improve the performance of OCN. 

Below we present two important ingredients when implementing our method to solve concrete problems, including those listed in Section \ref{exp}.

\subsection{Data sampling}
In this work, we assume the training data are collected from multiple trajectories of the dynamical system with randomly chosen initial points. To simulate this process, we generate the training data in our numerical experiments in the following way: 
\begin{itemize}
\item We first generate $m$ initial points from a specified distribution, say uniform distribution, over a domain in which we would like to learn the dynamical behavior of the solutions. Denote $y^{(j)}$ as the solution to (\ref{opcr}b) starting with the $j$-th initial point, the loss function in (\ref{opcr}) becomes
$$
J(\theta) = \sum_{j=1}^{m}\sum_{i=1}^{n}L_i(y^{(j)}(t_i)). 
$$
\item Starting with each initial point, we generate $\{x_i\}_{i=1}^n$ over time interval $[0, T]$ with $\Delta t=t_{i+1}-t_i$ for $i=1,...,n-1$ by solving the true dynamical system using a high-order ODE solver. For simplicity of notation, we assume the time interval $[0, T]$, the number of data points $n$, and the distance between two neighboring data points $\Delta t$ are the same for all trajectories.
\end{itemize}

\subsection{Batch training}
During training, each trajectory is divided into several mini-batches, and all batches of data are trained simultaneously. More precisely, for a trajectory data set of $\{x_i\}_{i=0}^{n}$, we divide it into $s$ batches: $\{x_{n_0},...,x_{n_1}\}$, ..., $\{x_{n_j},...,x_{n_{j+1}}\}$, ..., $\{x_{n_{s-1}},...,x_{n_{s}}\}$, where $n_0=0$ and $n_s=n$. 
From our experiments, we find that with fewer points in each batch, it takes less time to train the neural network to achieve a smaller training loss. Referring to Figure  \ref{fig:ocn}, the reason is that fewer points (or a shorter time interval) lead to less error accumulation due to the time discretization in step 2 and step 4, thus giving a more accurate gradient estimation in step 5. Hence, for a trajectory with $n+1$ points, we recommend dividing it into $n$ batches, with $2$ neighboring points in each batch.

\section{Error analysis}\label{ea}
In this section, we present theoretical results on the convergence behavior and error estimates for OCN. Note that the solution trajectory of (\ref{opcr}b) when an optimal parameter $\theta^*$ is obtained should be close to the solution trajectory of true dynamics (\ref{gf}).
Assume that $\nabla f(x)$ is Lipschitz continuous, and $x(t)$ is the unique solution to (\ref{gf}), and denote $y_k(t):= y(t; \theta_k)$ as the solution to (\ref{opcr}b) at the $k-$th iteration of training. These are functions evaluated at any point $t\in [0, T]=[t_0, t_n]$. We want to bound the error 
$$
e_k(t)=\|x(t)-y_k(t)\|.
$$
We will show that this error is bounded by the optimization error $J(\theta_K)$ and time step $O(\Delta t)$ with $\Delta t=\max_{0\leq i\leq n-1}|t_{i+1}-t_i|$.  

To quantify the errors and also control their propagation in time, we make the following assumptions: 

{\bf Assumption 1.} $f\in C^1(\Rb^d)$ and $\nabla f$ is Lipschitz continuous with constant $L_f$:
$$
\|\nabla f(x)-\nabla f(z)\|\leq L_f \|x-z\|, \quad \forall x, z\in \mathbb{R}^d. 
$$
Assumption 1 is a sufficient condition for the existence and uniqueness of the solution to (\ref{gf}). This is also used to control the truncation error in the discrete ODE (\ref{fe}).  

{\bf Assumption 2.}  $G\in C^1(\Rb^d\times \Rb^N)$ and there exist constant $L_{G_y}$ such that for any $\theta\in\mathcal{A}$,
$$
\|\partial_y G(y,\theta)-\partial_y G(z,\theta)\|\leq L_{G_y} \|y-z\|,\quad \forall y, z\in \mathbb{R}^d. 
$$
Assumption 2 plays a similar role for (\ref{opcr}b) as in Assumption 1 for the true dynamic system. Assumption 2 can be ensured by proper choices of activation functions in the construction of neural networks. In fact, we only need to take an activation function so that  $\sigma'$  is Lipschitz continuous. We note that the smoothness of the neural network may also be encouraged by the Lipschitz regularization \cite{LW+22}.


The main result is stated as follows. 

\begin{theorem}\label{thm2} Let Assumption 1 and 2 hold respectively on the regularity of $f$ and neural network $G$. 
Suppose that $\theta_k\in\mathcal{A}$ and $\mathcal{A}$ is bounded, where $\theta_k$ is generated using gradient descent (\ref{gd}) with gradient computed using Theorem \ref{thm1}, 
If $\Delta t=\max_{0\leq i\leq n-1} |t_{i+1}-t_i|\leq \frac{1}{2L_{G_y}}$, then 
\begin{align}\label{41}
    \max_{t \in [0, T]} \|x(t)-y_k(t)\|
    \leq C_1 (\sqrt{J(\theta_k)}+(\Delta t)^2).
\end{align}
 In addition, 
\begin{align}\label{42}
\max_i\|\nabla f(x_i)-\partial_y G(x_i, \theta_k)\|
\leq C_2 \left(\frac{\sqrt{J(\theta_k)}}{\Delta t}+\Delta t\right),
\end{align}
where $J(\theta_k)$ is the training loss defined by (\ref{loss}), $C_1, C_2$ are constants depending on the data, control set $\mathcal{A}$, and $L_f$ and $L_{G_y}$ in Assumptions 1 and 2. 
\end{theorem}
Due to space constraints, a detailed proof is relegated to Appendix \ref{pf-thm2}.  

Asymptotically, we expect $\lim_{k\to \infty}J(\theta_k)=J(\theta^*)$, which is zero or rather small, then the error in (\ref{41}) will ultimately be dominated by $(\Delta t)^2$, which is determined by how dense the data is collected over time. 

Without using any information on how dataset $\{x_i\}_{i=0}^{n}$ is sampled, the bound in (\ref{42}) may be the best possible one can get. However, if
the data is collected from solution trajectories of (\ref{gf}), then we expect $\frac{x_{i+1}-x_i}{\Delta t} \sim \dot x(t_i)$, which should be enforced to be close to $\dot y$ at $t_i$. With this consideration,  we may adopt an alternative loss function of form 
\begin{align}\label{opcr-v}
  \tilde J(\theta) = \sum_{i=1}^{n}\|y(t_i)-x_i\|^2 + \omega \sum_{i=1}^{n}\bigg\|\frac{x_{i}-x_{i-1}}{\Delta t}+\partial_y G(y(t_{i-1}),\theta)\bigg\|^2,
\end{align} 
where $\omega>0$ is a weighting parameter. 

\begin{theorem}\label{thm2-v} Under the same conditions as in Theorem \ref{thm2}, with loss function (\ref{opcr-v}) used in training, the error bound (\ref{41}) still holds, and  
\begin{align}\label{42-v}
\max_i\|\nabla f(x_i)-\partial_y G(x_i, \theta_k)\|
\leq C_2 (\sqrt{\tilde J(\theta_k)}+\Delta t),
\end{align}
where $\tilde J(\theta_k)$ is the training loss defined by (\ref{opcr-v}), $C_2$ are constants depending on the observed data, control set $\mathcal{A}$, and  $L_f$ and $L_{G_y}$ in Assumptions 1 and 2. 
\end{theorem}
The proof of this theorem is similar,  we defer details to  Appendix \ref{pf-thm2-v}.

\section{Time-discretization}\label{discre}
In this section, we discuss how to discretize system (\ref{00}) in order to accurately evaluate the gradient (\ref{grad-}).  
One approach is to integrate an augmented system backward in time, as in the original implementation of the neural ODEs \cite{CR+18}.  However, there are some observed drawbacks: possible instability in solving (\ref{opcr}b) backward in time; the computational cost is twice more than the ordinary backpropagation algorithm; numerical errors can also harm the accuracy of the gradient estimation. 

\subsection{Symplectic integrator} In order to enhance the accuracy of the gradient estimation with (\ref{grad-}), we seek a time-discretization that can conserve some time-invariants. In system (\ref{00}), one can verify that there are two time-invariants in each interval $t\in (t_i, t_{i+1}]$, 
\begin{align*}
& H =-\partial_y G(y,\theta)p, \\
& S =\delta^\top p, \quad \delta(t):= \frac{\partial y(t)}{\partial y(0)}. 
\end{align*}
Here $y(0)$ serves as the initial data for the forward problem, and $y$ is the corresponding flow map $y=\phi(t; y(0))$. The first quantity $H$ is a Hamiltonian. Typically, one can only hope to conserve certain modified Hamiltonian by a high-order ODE solver. The second quantity $S$ is bilinear and associated with the symplectic structure of the coupled system (\ref{00}). In fact, by the chain rule, we have 
\begin{equation}\label{Sp}
\frac{d}{dt}S  =\dot \delta^\top  p
+\delta^\top \dot p
 = (-\partial_y^2 G(y;\theta)\delta)^\top p +
\delta^\top (-\partial_y^2 G(y, \theta))^\top p)=0.    
\end{equation} 
As shown in \cite{MM21}, a partitioned Runge-Kutta method can be formulated  to conserve $S$ at the discrete level.  

To be more concrete, we discretize the forward equation by a Runge-Kutta (RK) method. Let $t_l$, $\tau_l=t_{l+1}-t_l$, $y_l$ denote the $l$-th time step, step size, and state, respectively. RK method with $s$ stages has the following  form
\begin{equation}\label{arkf}
\begin{aligned}
y_{l+1} &= y_l + \tau_l \sum_{i=1}^{s}b_i g_{li}, \\
g_{li} &:= -\partial_y G(y_{li},\theta), \\
y_{li} &= y_l + \tau_l \sum_{j=1}^{s}a_{ij}g_{lj},
\end{aligned}    
\end{equation}
where $a_{ij}, b_i$ are the RK coefficients. In the case $b_i\not=0$ for all $i\in\{1\cdots s\}$, the backward problem is solved by another RK method with the same step size as that used for the system state $y$, with RK coefficients: $A_{ij}$ and $B_i$. Such a partitioned RK method for system (\ref{00}) can be shown to conserve $S$ as long as 
$$
b_iA_{ij}+B_ia_{ji}-b_iB_j=0 \; \text{for} \; i, j=1,\cdots, s, \; \text{and} \; B_i=b_i\; \text{for}\; i=1,\cdots s.
$$
For RK methods with some $b_i=0$, a modified scheme for the backward problem can be formulated as 
\begin{equation}\label{arkb}
\begin{aligned}
p_{l} &= p_{l+1} - \tau_l \sum_{i=1}^{s}\tilde b_i h_{li},\\
h_{li} &:= \partial^2_y G(y_{li},\theta)^\top p_{li},\\
p_{li} &: = 
\begin{cases}
p_{l+1} - \tau_l \sum_{j=1}^{s}\tilde b_j\frac{a_{ji}}{b_i}h_{lj}, & \text{if}\quad b_i \neq 0,\\
- \sum_{j=1}^{s}\tilde b_ja_{ji}h_{lj}, & \text{if}\quad b_i =  0,
\end{cases}
\end{aligned}    
\end{equation}
where $\tilde b_i=b_i$ if $b_i\neq 0$ else $\tilde b_i=\tau_l$. Note that  (\ref{arkb}) is explicit backward in time as long as the RK method in (\ref{arkf}) is explicit forward in time, which is the case when $a_{ij}=0$ for $j\geq i$.

\begin{theorem}\label{Sinv}
If the forward problem (\ref{00}a) is solved by (\ref{arkf}), and for each time interval $(t_{i-1},t_{i}]$, where $i=n,...,1$, the backward problem (\ref{00}b) is solved by (\ref{arkb}), then in each time interval $(t_{i-1},t_{i}]$, the quantity $\delta ^\top p$ is conversed, i.e., $\delta_{l+1}^\top p_{l+1}=\delta_{l}^\top p_{l}$ for all $l\geq 0$, where $\delta_{l+1}=\frac{\partial y_{l+1}}{\partial y(0)}$.
\end{theorem}
The proof is deferred to Appendix \ref{pf-Sinv}.

In our experiments, we use Dopri5 (5th-order Dormand–Prince method) \cite{DP80}, an RK method with adaptive step size, to discretize the forward problem. It takes the form of (\ref{arkf}) with $s=7$, and $b_2=b_7=0$. The backward problem is discretized using  (\ref{arkb}). 

\section{Experimental results}\label{exp}
In this section, we test the proposed method on several canonical systems.\footnote{The code is available at \url{https://github.com/txping/OCN}.} For all experiments, we use feed-forward neural networks with the tanh activation function. The detailed structure of the neural network applied for each problem is provided in corresponding subsections. All the weights are initialized randomly from Gaussian distributions, and all the biases are initialized to zero.

After the neural network is well trained, we generate $\{y(t_i)\}_{i=1}^n$ from the learned dynamics $\dot y=-\partial_y G(y,\cdot)$ (or $\dot y= G(y,\cdot)$) and compare it against the observed data $\{x_i\}_{i=1}^n$. For the first three examples, the comparison between $G(x_i,\cdot)$ and $f(x_i)$ is given. For experiments on the gradient flow problem, we also verify the generalization performance of OCN by applying it to testing data, which are some initial points generated randomly over the same domain and do not appear in the dataset used for training. 

For each experiment, we provide the true dynamical system, which is used to generate the observed data and verify the performance of the trained models, but in no way facilitates the neural network approximation.

\subsection{Linear gradient flow} 
For this example, the observed data is collected on solution trajectories to
\begin{align*}
\dot x_1 = -2x_1 - x_2,\\
\dot x_2 = - x_1 - 2x_2.
\end{align*}
This is of form $\dot x=-\nabla f(x)$ with
\begin{equation}\label{lgf-f}
f(x_1, x_2) = x^2_1 + x_1x_2 + x^2_2.  
\end{equation}
This system has critical point $(0,0)$ as a stable node. All solution trajectories tend to $(0,0)$ as $t\to\infty$. We want to extract $f$ from the training data, which is sampled from $8$ trajectories on domain $[-2,2]\times[-2,2]$ with time interval $[0,5]$ and time step $\Delta t=0.05$. The neural network $G$ used to approximate $f$ in (\ref{lgf-f}) has $2$ hidden layers of $50$ neurons.

\begin{figure}[ht]
\begin{subfigure}[b]{0.33\linewidth}
\centering
\includegraphics[width=1\linewidth]{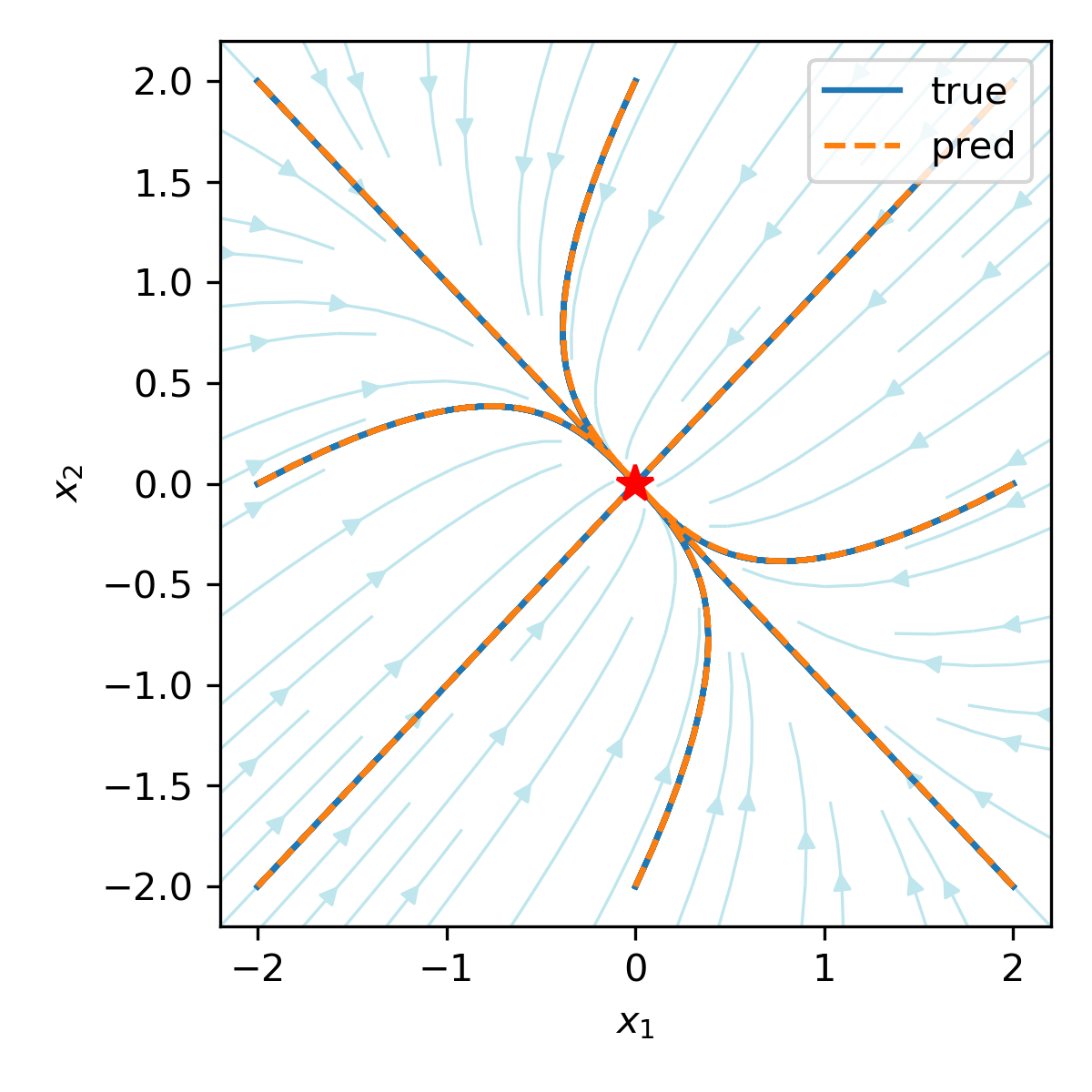}
\caption{Training result}
\end{subfigure}%
\begin{subfigure}[b]{0.33\linewidth}
\centering
\includegraphics[width=1\linewidth]{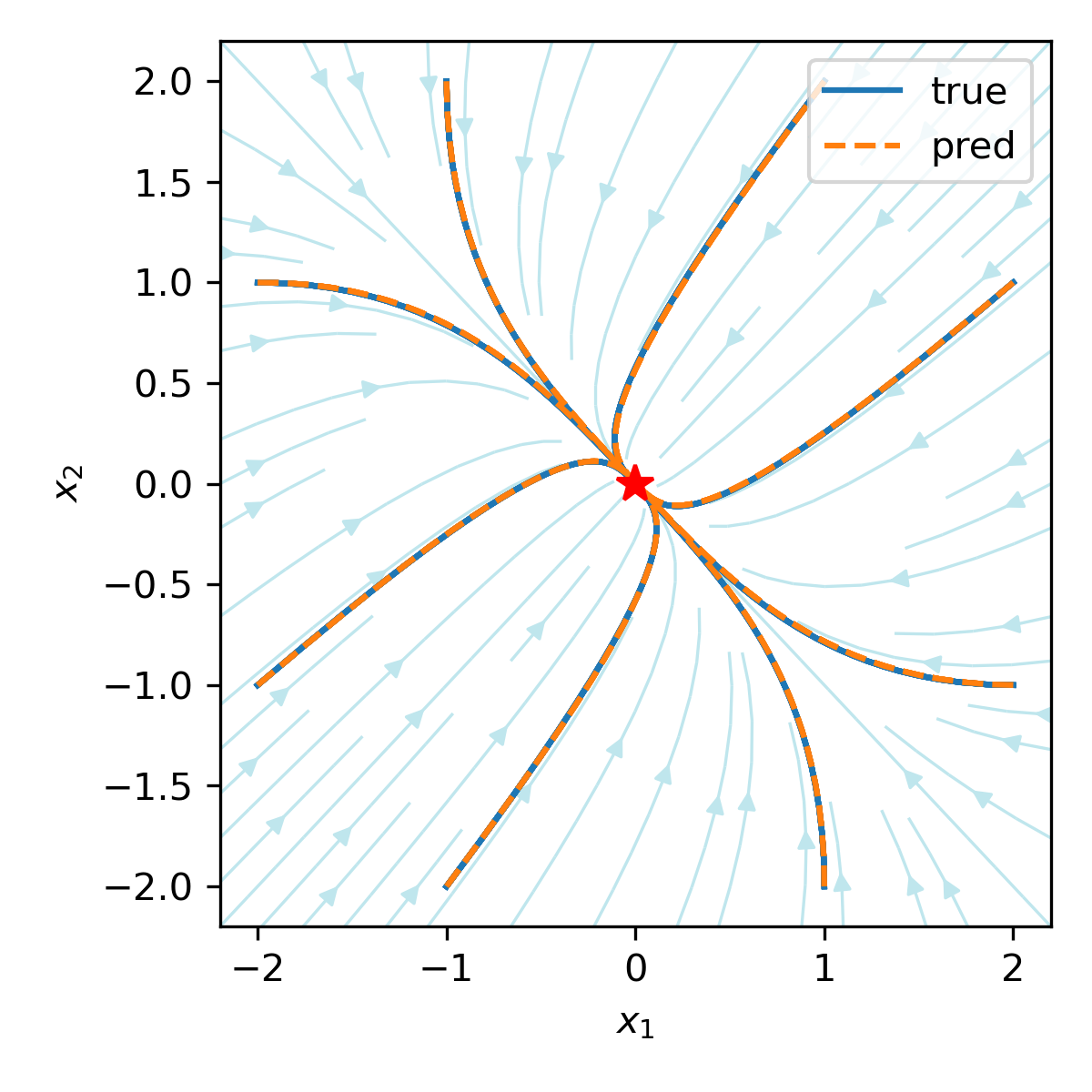}
\caption{Test result}
\end{subfigure}%
\begin{subfigure}[b]{0.33\linewidth}
\centering
\includegraphics[width=1\linewidth]{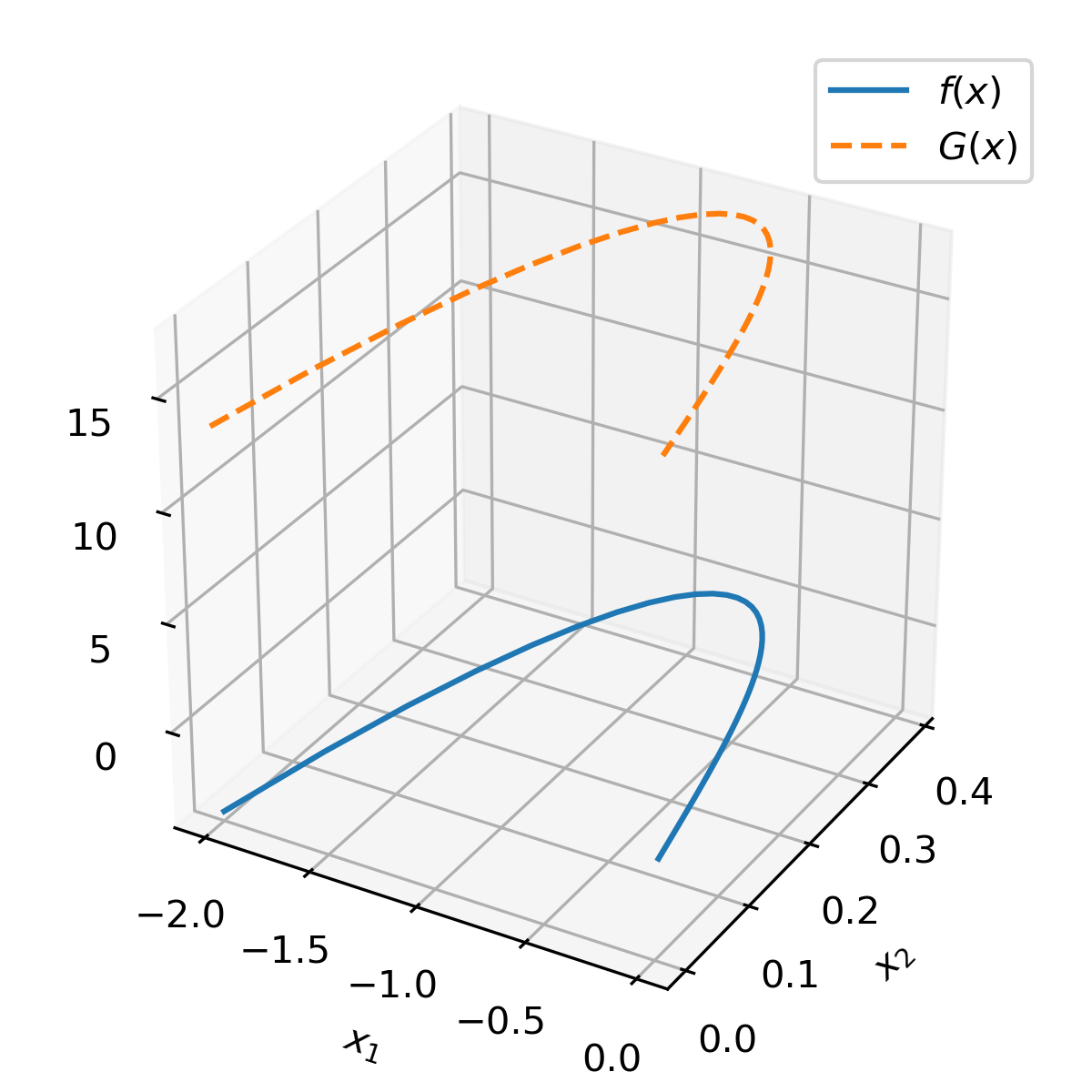}
\caption{Function profile}
\end{subfigure}%
\caption{Results of the linear gradient flow. For (a) and (b), the star represents the minimizer of $f$ in (\ref{lgf-f}). 
}
\label{fig:eg1}
\end{figure}

The training and testing results are presented in Figure \ref{fig:eg1} (a) and (b), respectively. It can be seen that all trajectories generated by OCN match the observed data generated by the true dynamical system well. 

Figure \ref{fig:eg1} (c) is a comparison between the true governing function $f(x)$ and the trained neural network $G(x,\cdot)$, where $x$ represents the training data set $\{x_i\}$. 
$G(x,\cdot)$ is an affine translation of the true function because the original problem (\ref{gf}) is uniquely determined up to a constant, $f+c$ for any constant $c$. For $G(x,\cdot)$ that satisfies (\ref{gf}), $G(x,\cdot)+c$ also satisfies (\ref{gf}) for any constant $c$. 

\subsection{Nonlinear gradient flow}
For this example, the observed data is collected on solution trajectories to
\begin{equation}\label{ngf}
\begin{aligned}
\dot x_1 &= -\cos(x_1)\cos(x_2),\\
\dot x_2 &= \sin(x_1)\sin(x_2).
\end{aligned}
\end{equation}
This is of form $\dot x= -\nabla f(x)$ with
\begin{equation}\label{ngf-f}
f(x_1, x_2) = \sin(x_1)\cos(x_2). 
\end{equation}
This system has three types of nodes -- stable nodes, unstable nodes, and saddle points -- spread over the domain in a staggered pattern. Stable nodes at $[(k_1+\frac{1}{2})\pi, k_2\pi]$ where $k_1$ and $k_2$ have opposite parity; unstable nodes at $[(k_1+\frac{1}{2})\pi, k_2\pi]$ where $k_1$ and $k_2$ have the same parity; saddle points at $[k_3\pi, (k_4+\frac{1}{2})\pi]$. The training data consists of $24$ trajectories sampled from domain $[-6,6]\times[-4,6]$ with time interval $[0,8]$ and $\Delta t=0.05$. The neural network $G$ used to approximate $f$ in (\ref{ngf-f}) has $2$ hidden layers of $200$ neurons.

\begin{figure}[ht]
\begin{subfigure}[b]{0.33\linewidth}
\centering
\includegraphics[width=1\linewidth]{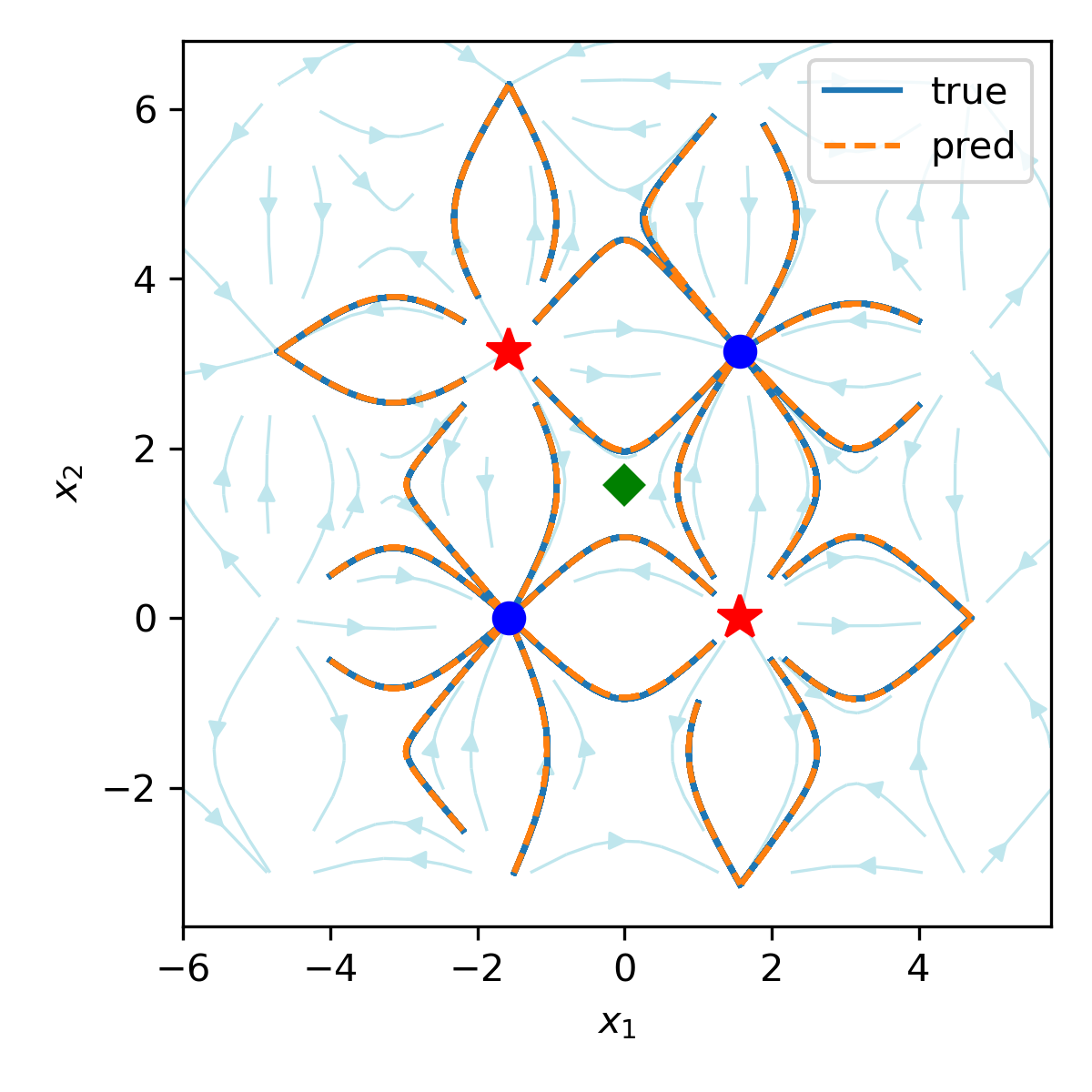}
\caption{Training result}
\end{subfigure}%
\begin{subfigure}[b]{0.33\linewidth}
\centering
\includegraphics[width=1\linewidth]{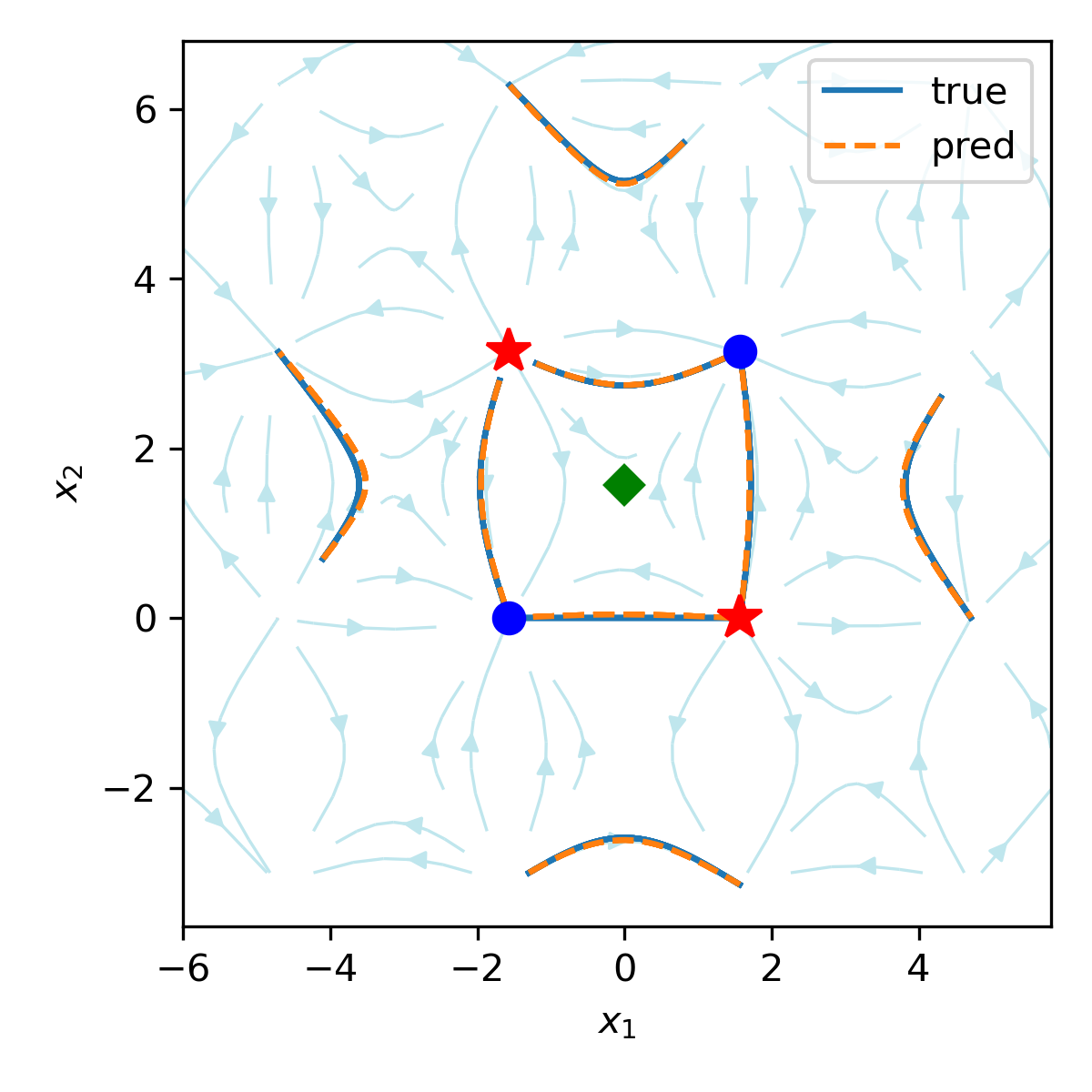}
\caption{Test result}
\end{subfigure}%
\begin{subfigure}[b]{0.33\linewidth}
\centering
\includegraphics[width=1\linewidth]{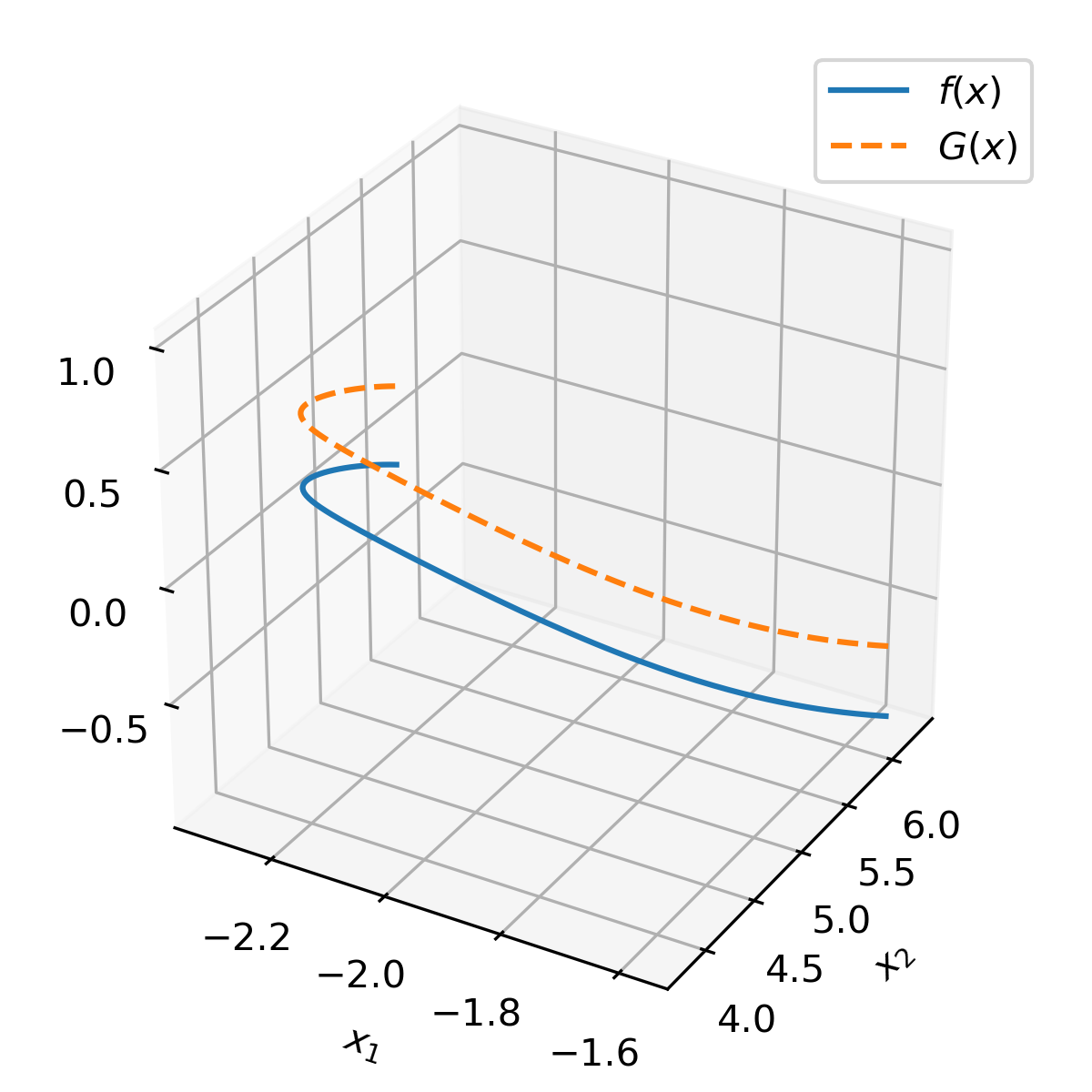}
\caption{Function profile}
\end{subfigure}%
\caption{Results of the nonlinear gradient flow. The stars represent unstable nodes, the circles represent stable nodes, and the squares represent saddle points.}
\label{fig:eg2}
\end{figure}

The training results are presented in Figure \ref{fig:eg2} (a). We observe that for trajectories around different types of nodes, either diverging from sources or converging to sinks, OCN fits the training data well. 

The performance of OCN on test data is shown in Figure \ref{fig:eg2} (b). The test data is composed of $8$ initial points, among which $4$ initial points (in the center of the figure) correspond to trajectories that have a similar pattern to that of the training data; another $4$ initial points correspond to trajectories whose dynamic behavior is different from that of the training data. For both types of initial points, OCN recovers the true trajectories well.

\subsection{Damped pendulum}\label{nodetest} To illustrate that our method is well applicable to general ODE systems, we consider the pendulum problem, which has the form of $\dot x(t)=F(x(t))$. Specifically,
\begin{align*}
&\dot x_1 = x_2,\\
&\dot x_2 = -0.2 x_2-8.91 \sin(x_1).
\end{align*}
Here $x_1$ is the angular displacement, and $x_2$ is the angular velocity. This is a damped system that obeys a dissipation law: 
$$
\frac{d}{dt}\bigg(\frac{x^2_2}{2} +8.91(1-\cos(x_1))\bigg) = -0.2 x^2_2 \leq 0.
$$
The critical point $(0,0)$ is a stable focus. The training data is collected from $1$ trajectory starting from $[-1,-1]$ within time interval $[0,5]$ and time step $\Delta t=0.05$. The neural network $G$ used to approximate $f$ has $1$ hidden layer of $100$ neurons.

\begin{figure}[ht]
\begin{subfigure}[b]{0.33\linewidth}
\centering
\includegraphics[width=1\linewidth]{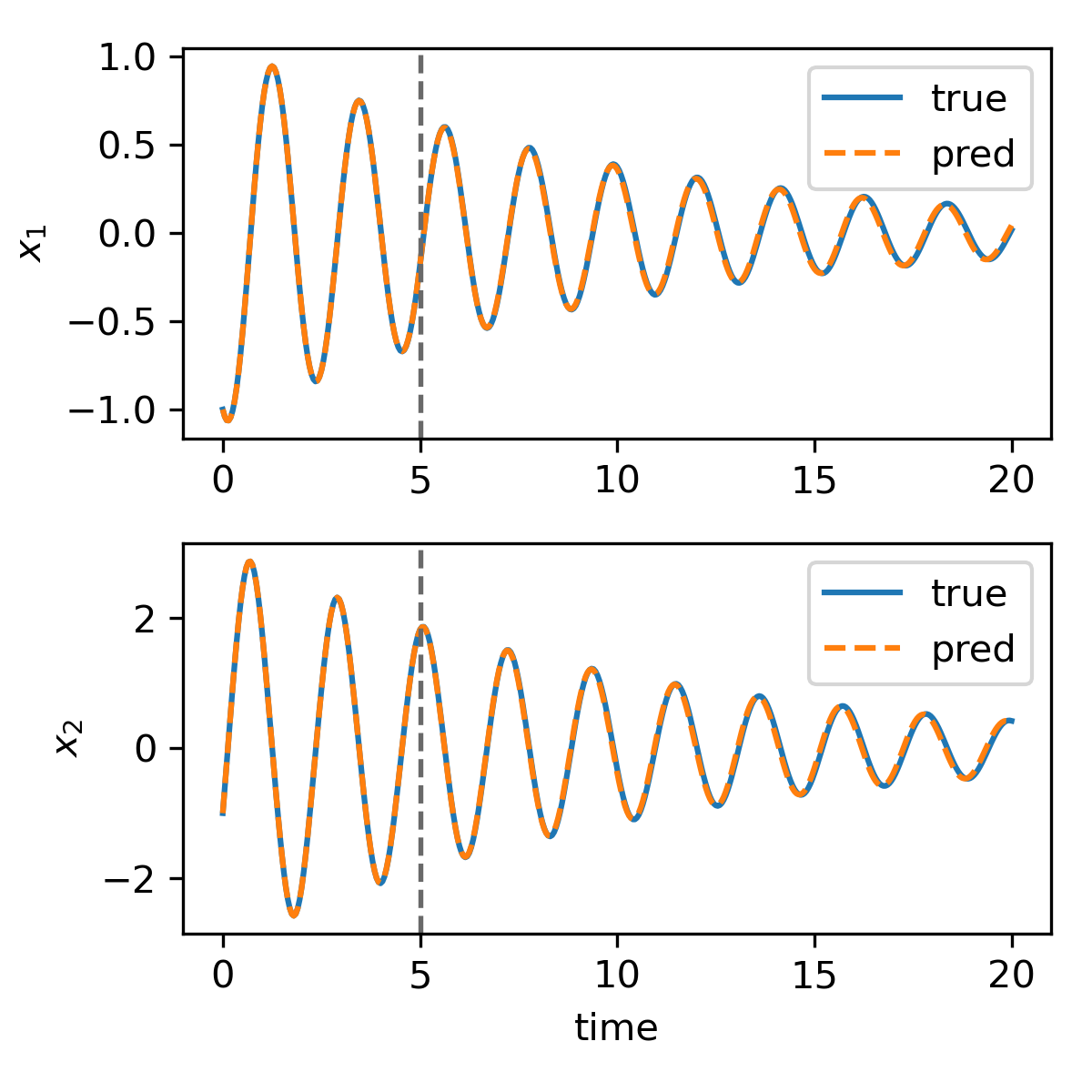}
\caption{Trajectory}
\end{subfigure}%
\begin{subfigure}[b]{0.33\linewidth}
\centering
\includegraphics[width=1\linewidth]{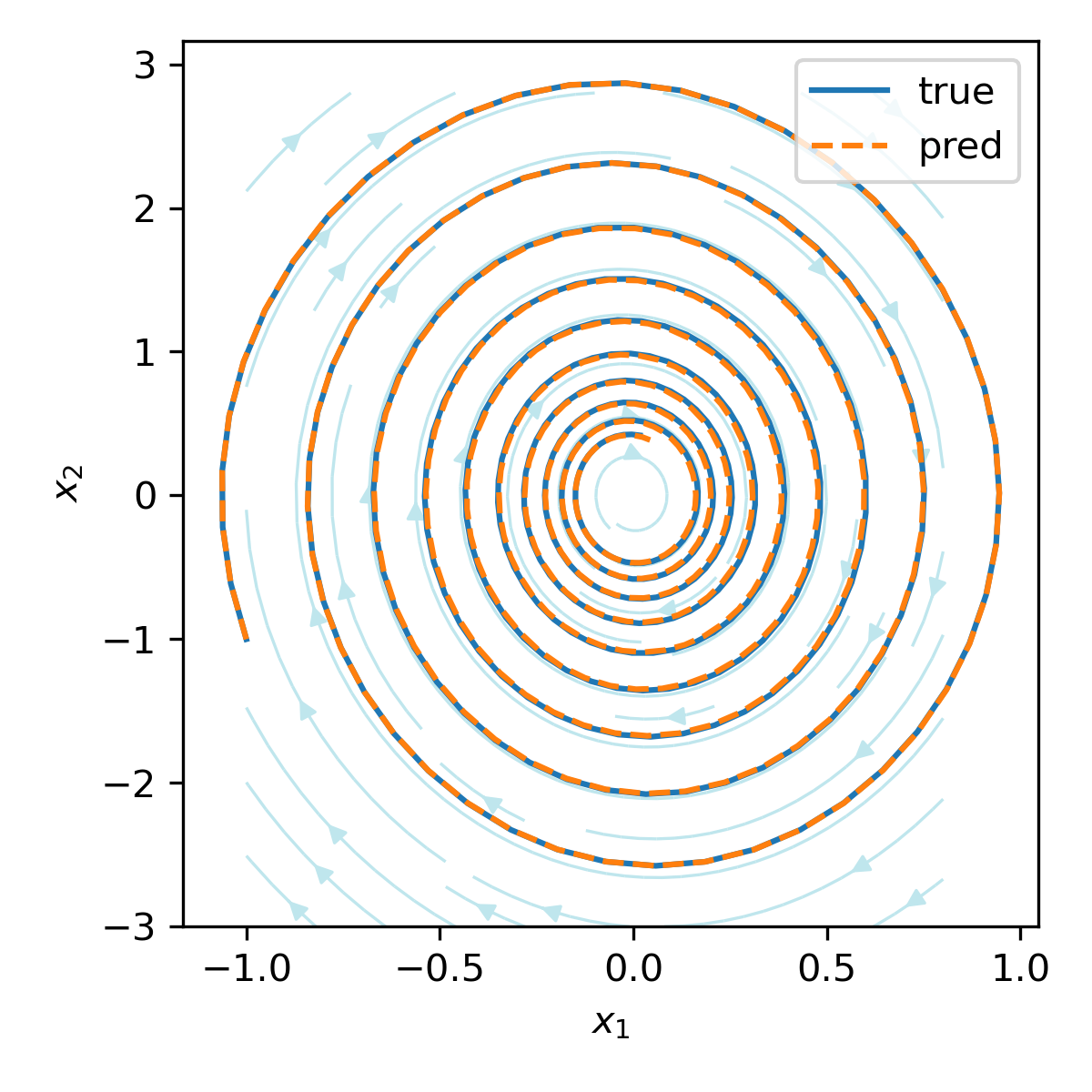}
\caption{Phase portrait}
\end{subfigure}%
\begin{subfigure}[b]{0.33\linewidth}
\centering
\includegraphics[width=1\linewidth]{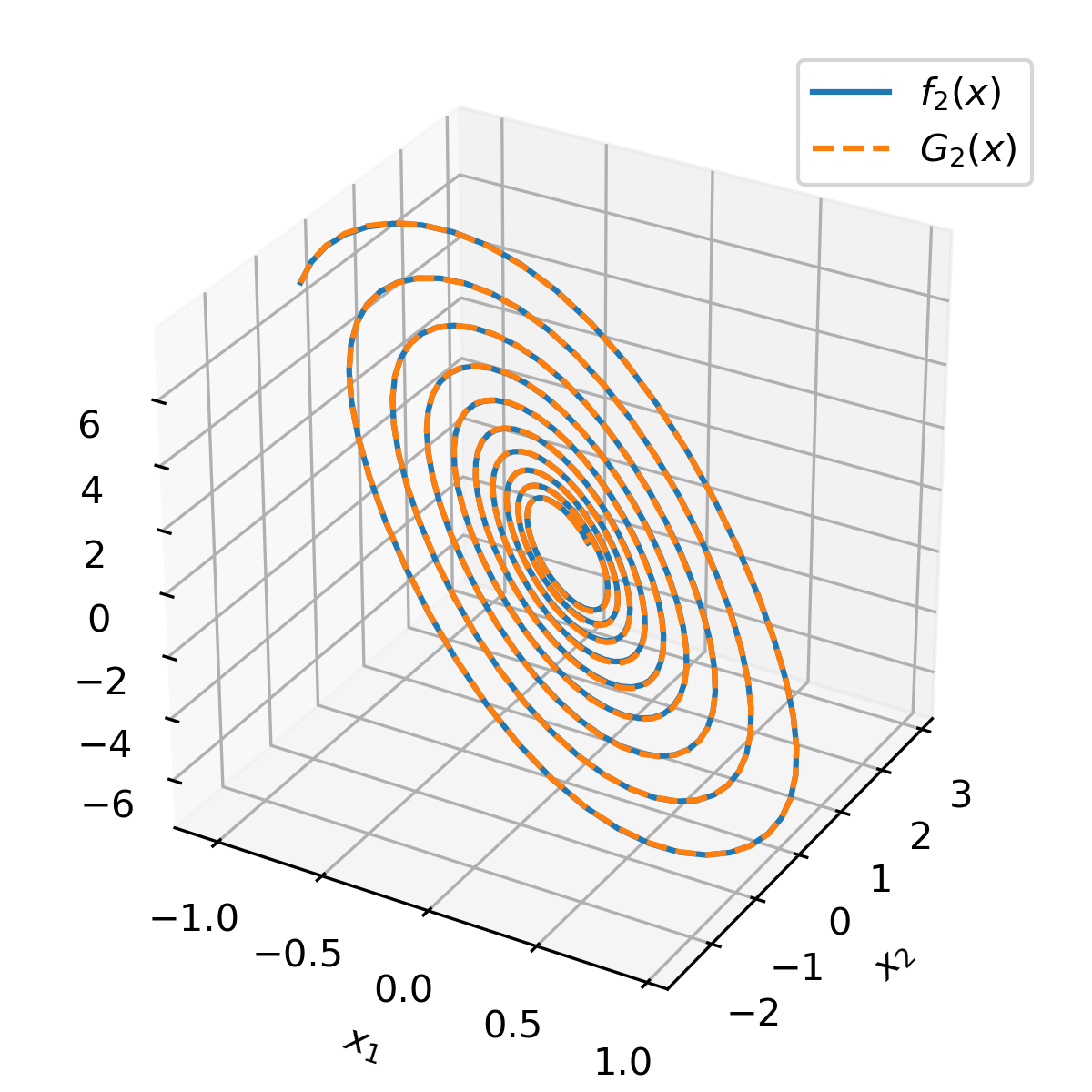}
\caption{Function profile}
\end{subfigure}%
\caption{Results of the nonlinear ODE system. The results on $[0, 5]$ show the performance of OCN on the training data; the results on $[5, 20]$ show the prediction performance of OCN. For (c), $f_2(x)=-0.2 x_2-8.91 \sin(x_1)$.}
\label{fig:eg3}
\end{figure}

After finishing training, we generate a trajectory over $[0,20]$ to examine the relatively long-term prediction performance of OCN. The results are presented in Figure \ref{fig:eg3}. We observe accurate fitting between the true trajectory and the trajectory generated by OCN, even on a time interval that is much longer than what is used for training.


\subsection{Lorenz system}
We demonstrate our method on the nonlinear Lorenz system \cite{Lor63}:
\begin{equation}\label{Lorenz}
\begin{aligned}
&\dot x = \sigma(y-x),\\
&\dot y = x(\rho - z) - y,\\
&\dot z = xy - \beta z.
\end{aligned}    
\end{equation}
The dynamics are very rich for different choices of parameters $(\sigma,\rho,\beta)$. The well-known Lorenz attractor shows up for $(\sigma,\rho,\beta)=(10, 28, 8/3)$. 
For this example, the neural network $G$ used to approximate $f$ has $3$ hidden layers of $300$ neurons. The detailed experimental setup is given below; see also Table \ref{tb} for a summary of the results. 

\subsubsection{Generalization performance} 
\begin{figure}[ht]
\begin{subfigure}[b]{0.25\linewidth}
\centering
\includegraphics[width=1\linewidth]{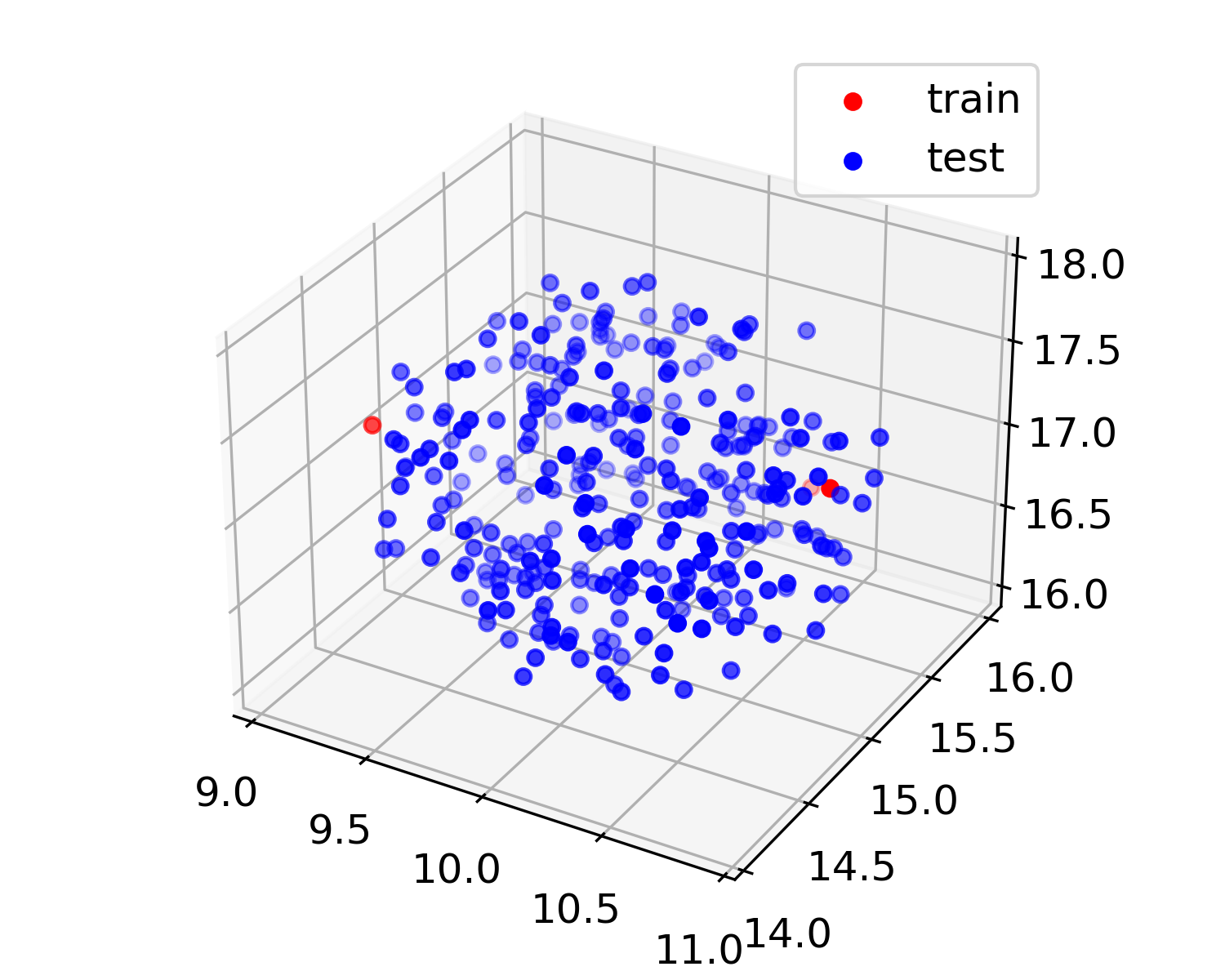}
\caption{Initial points}
\end{subfigure}%
\begin{subfigure}[b]{0.25\linewidth}
\centering
\includegraphics[width=1\linewidth]{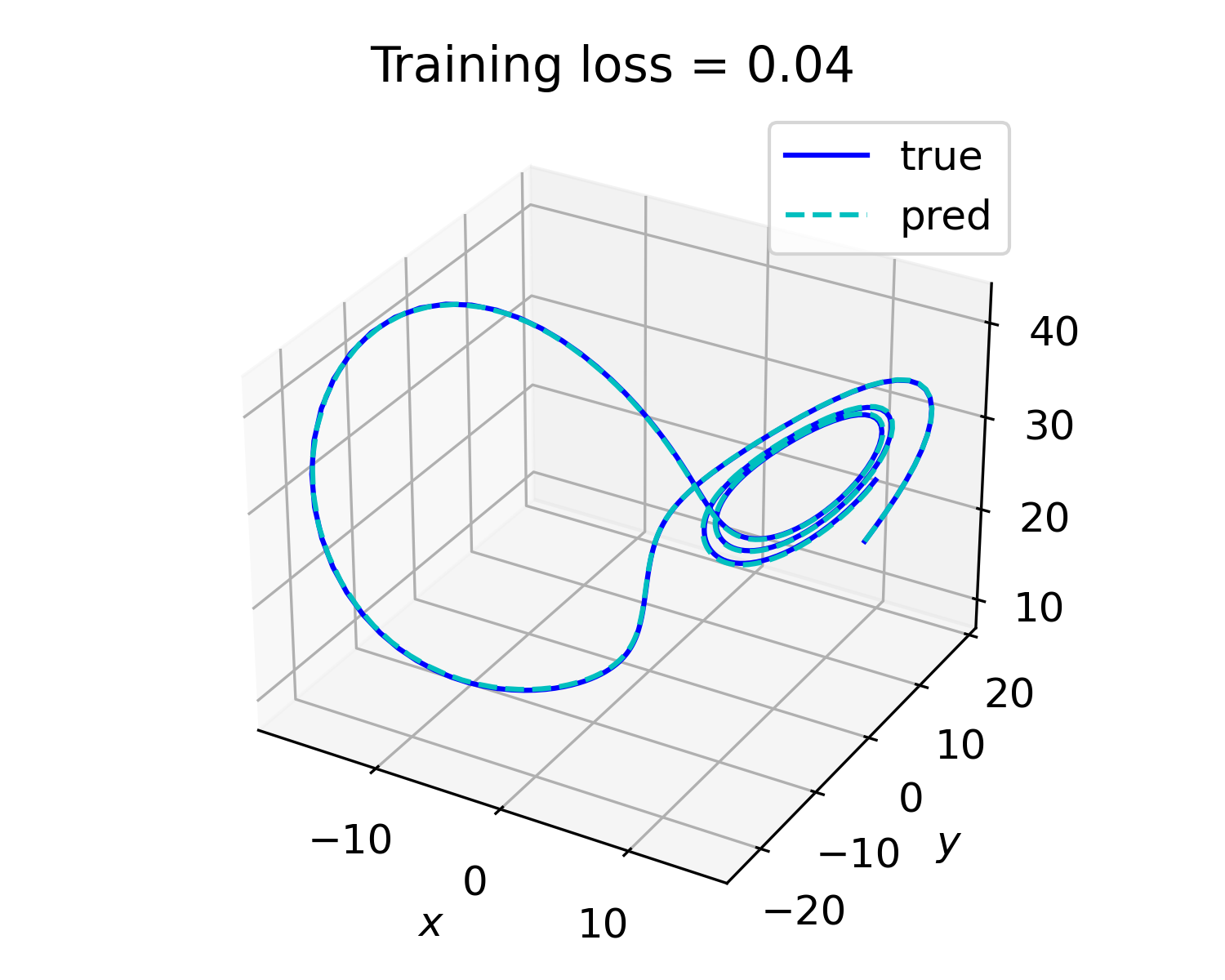}
\caption{Train}
\end{subfigure}%
\begin{subfigure}[b]{0.25\linewidth}
\centering
\includegraphics[width=1\linewidth]{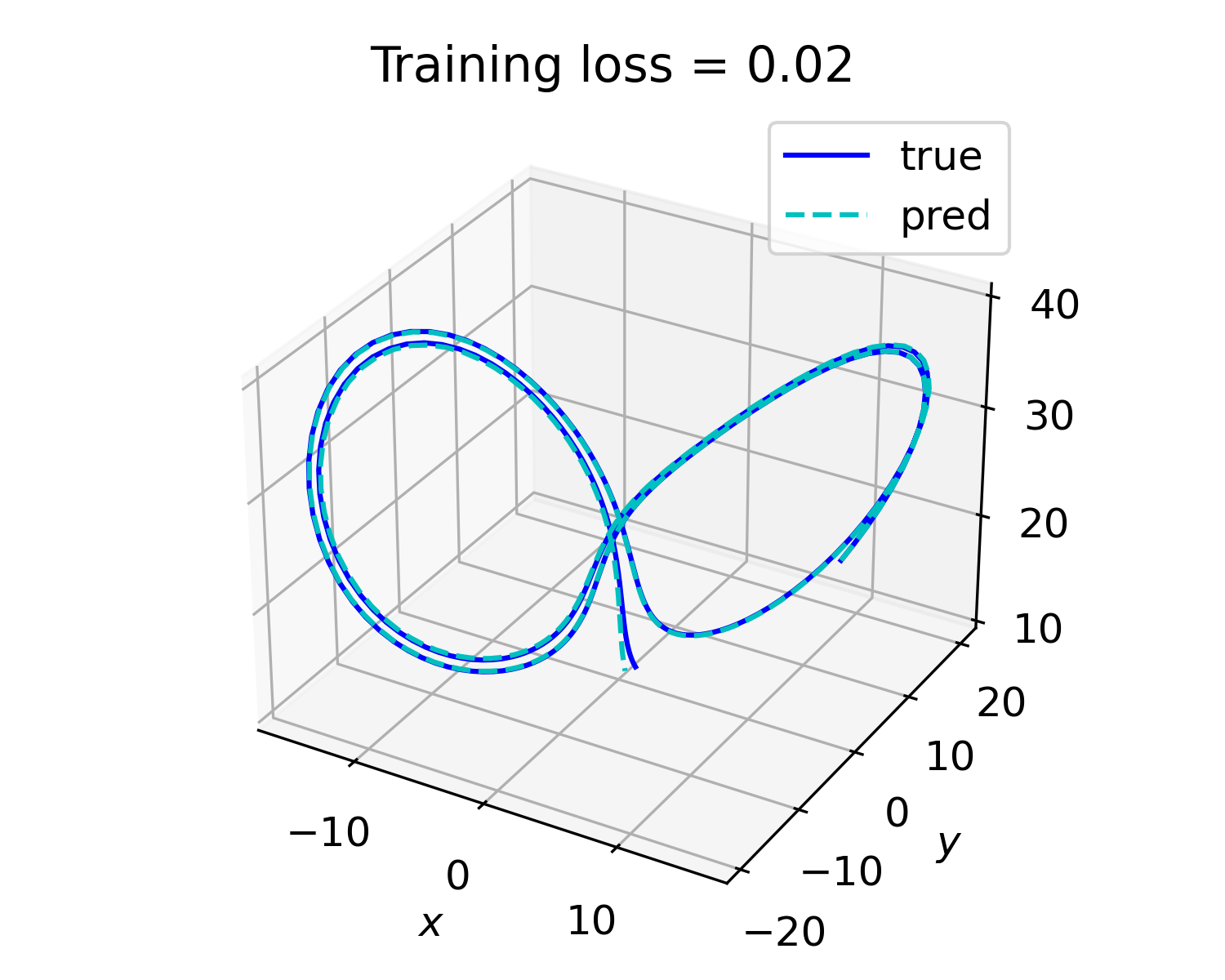}
\caption{Train}
\end{subfigure}%
\begin{subfigure}[b]{0.25\linewidth}
\centering
\includegraphics[width=1\linewidth]{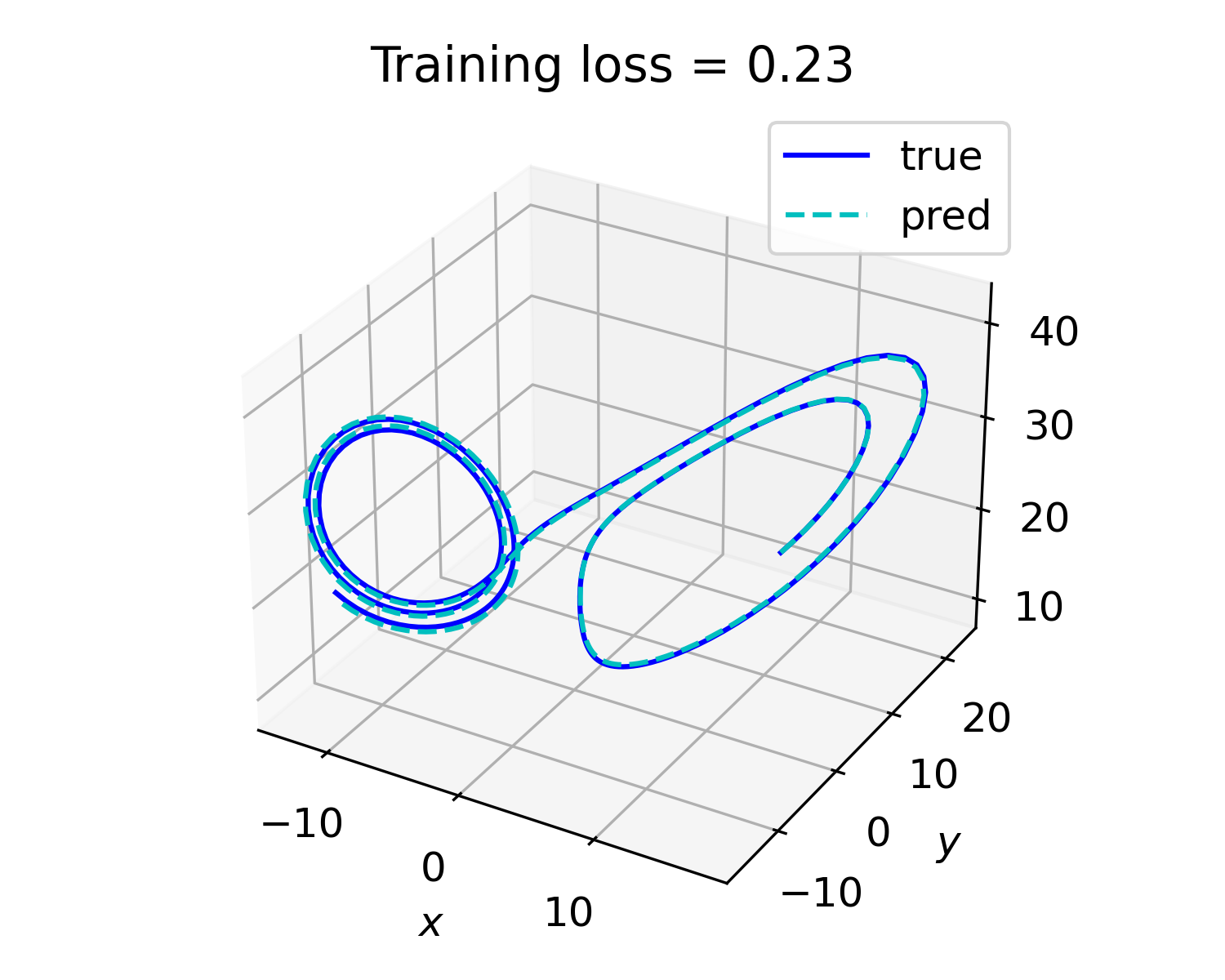}
\caption{Train}
\end{subfigure}%
\newline 
\begin{subfigure}[b]{0.24\linewidth}
\centering
\includegraphics[width=1\linewidth]{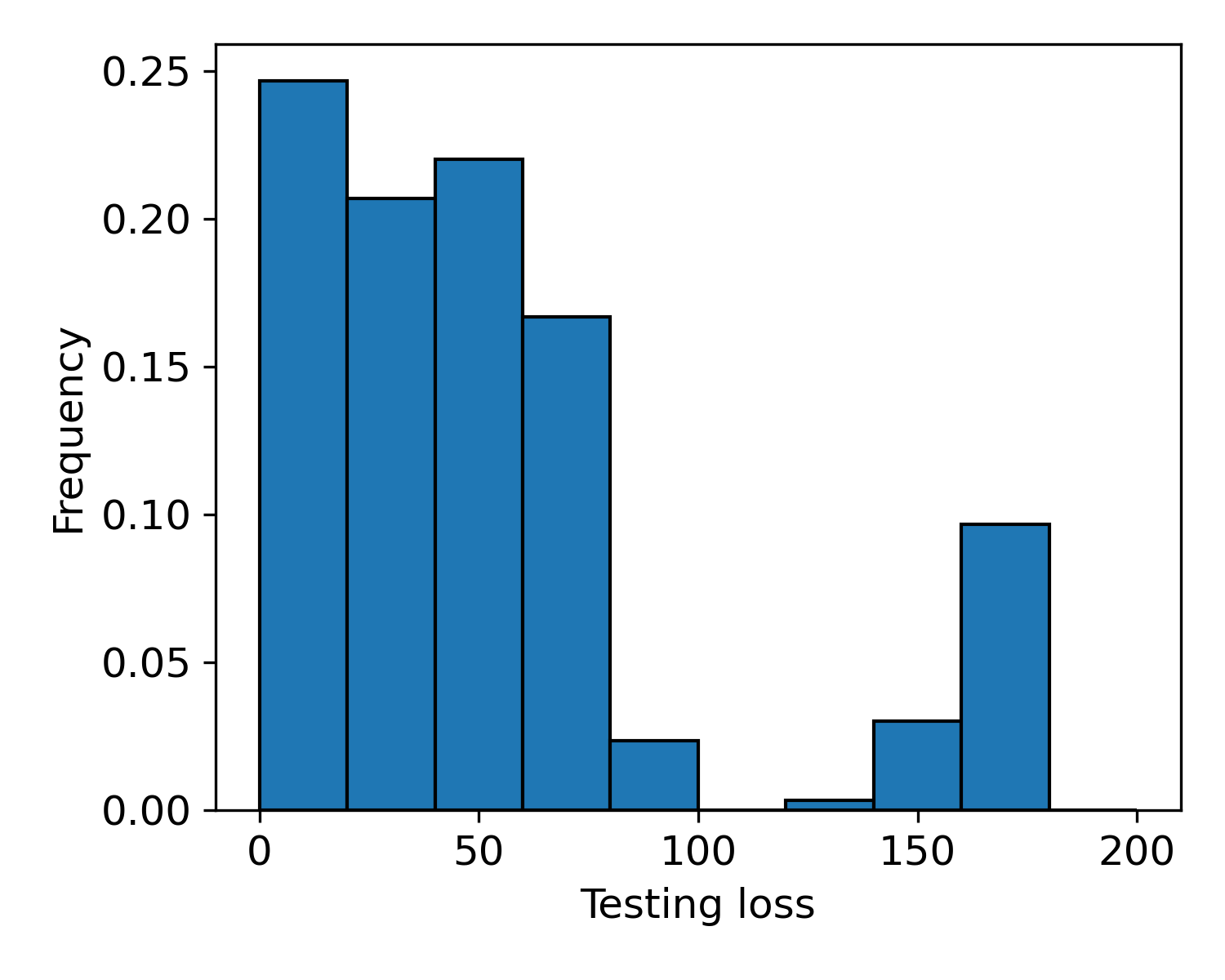}
\caption{Histogram}
\end{subfigure}%
\begin{subfigure}[b]{0.25\linewidth}
\centering
\includegraphics[width=1\linewidth]{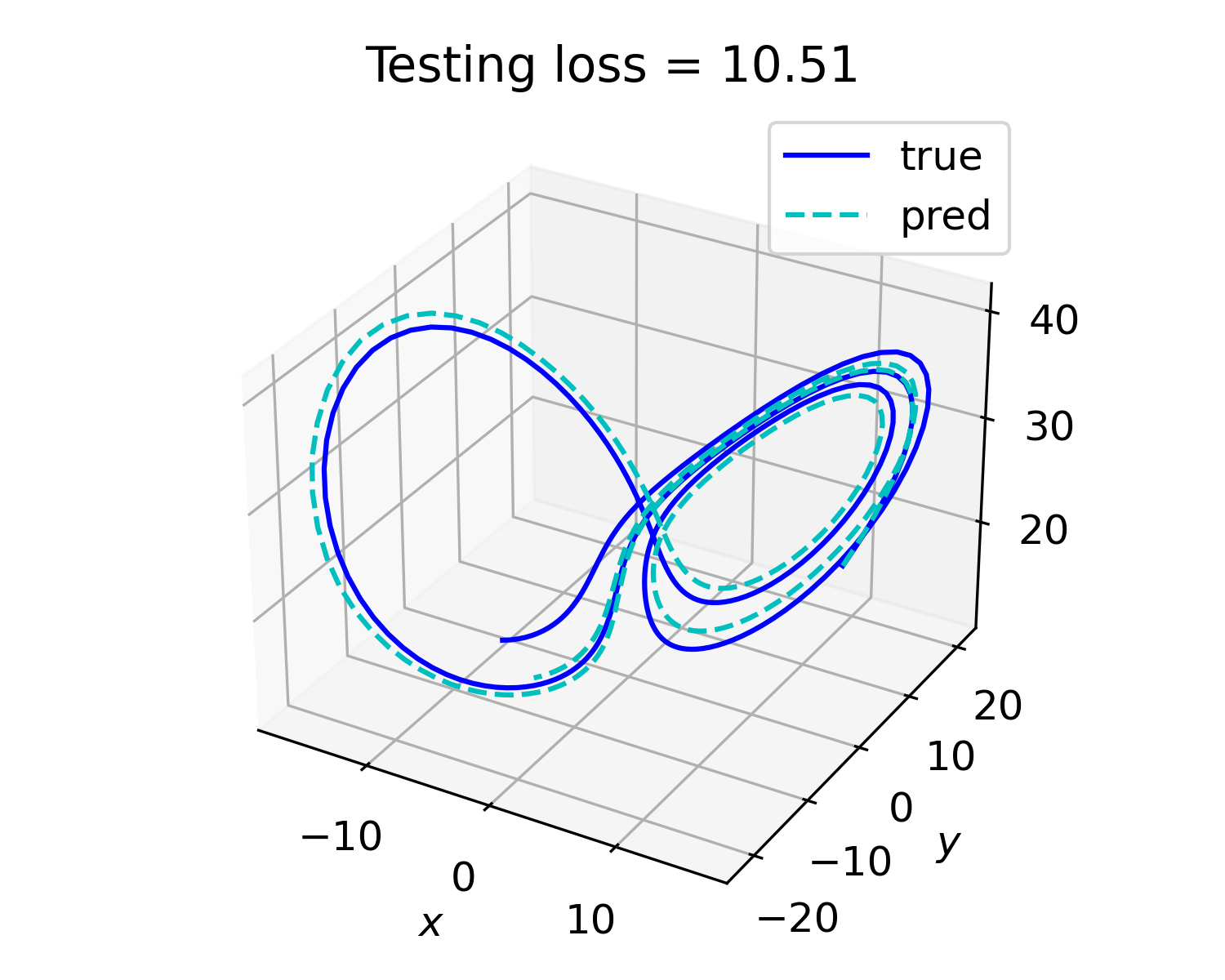}
\caption{Test (unseen)}
\end{subfigure}%
\begin{subfigure}[b]{0.25\linewidth}
\centering
\includegraphics[width=1\linewidth]{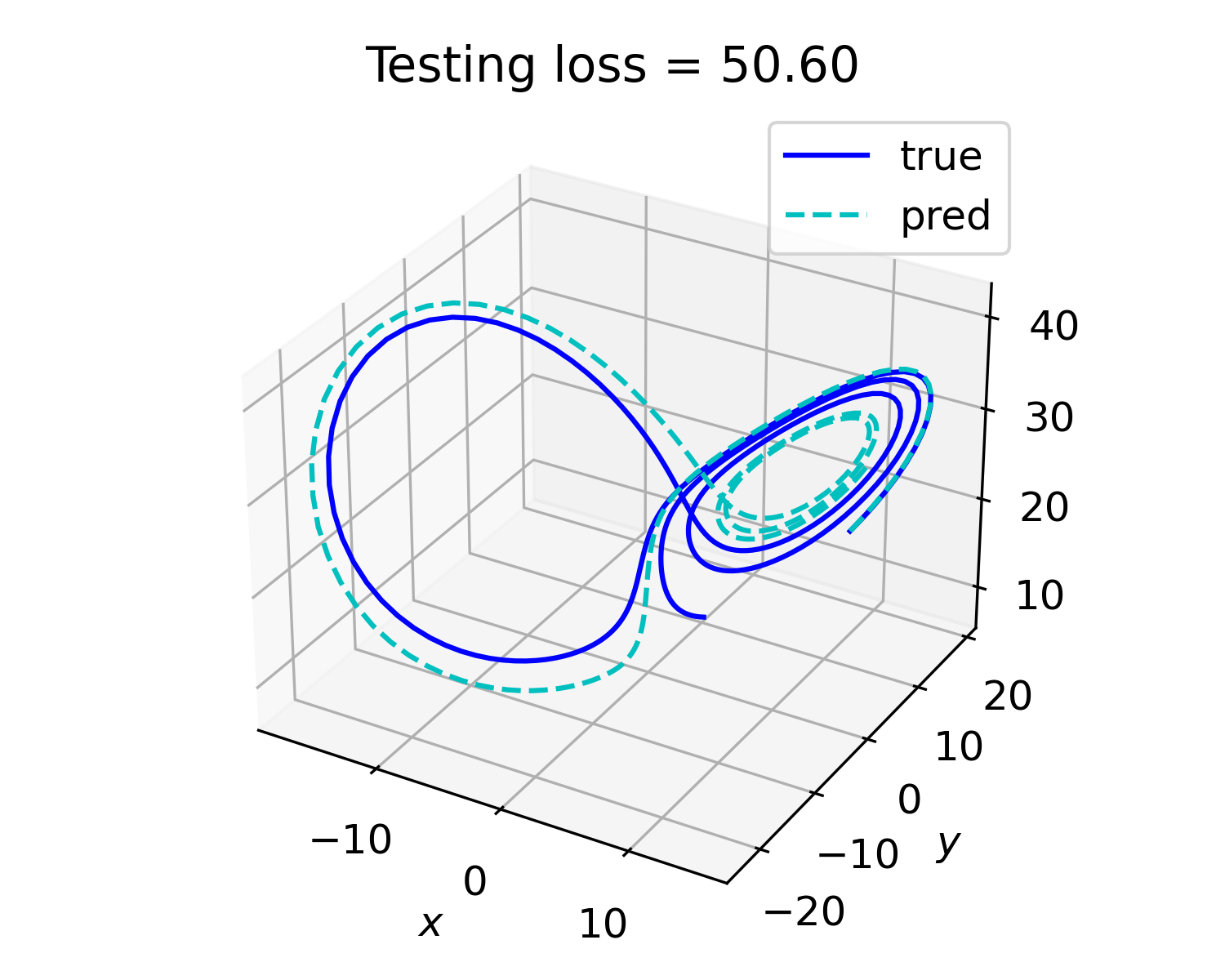}
\caption{Test (unseen)}
\end{subfigure}
\begin{subfigure}[b]{0.25\linewidth}
\centering
\includegraphics[width=1\linewidth]{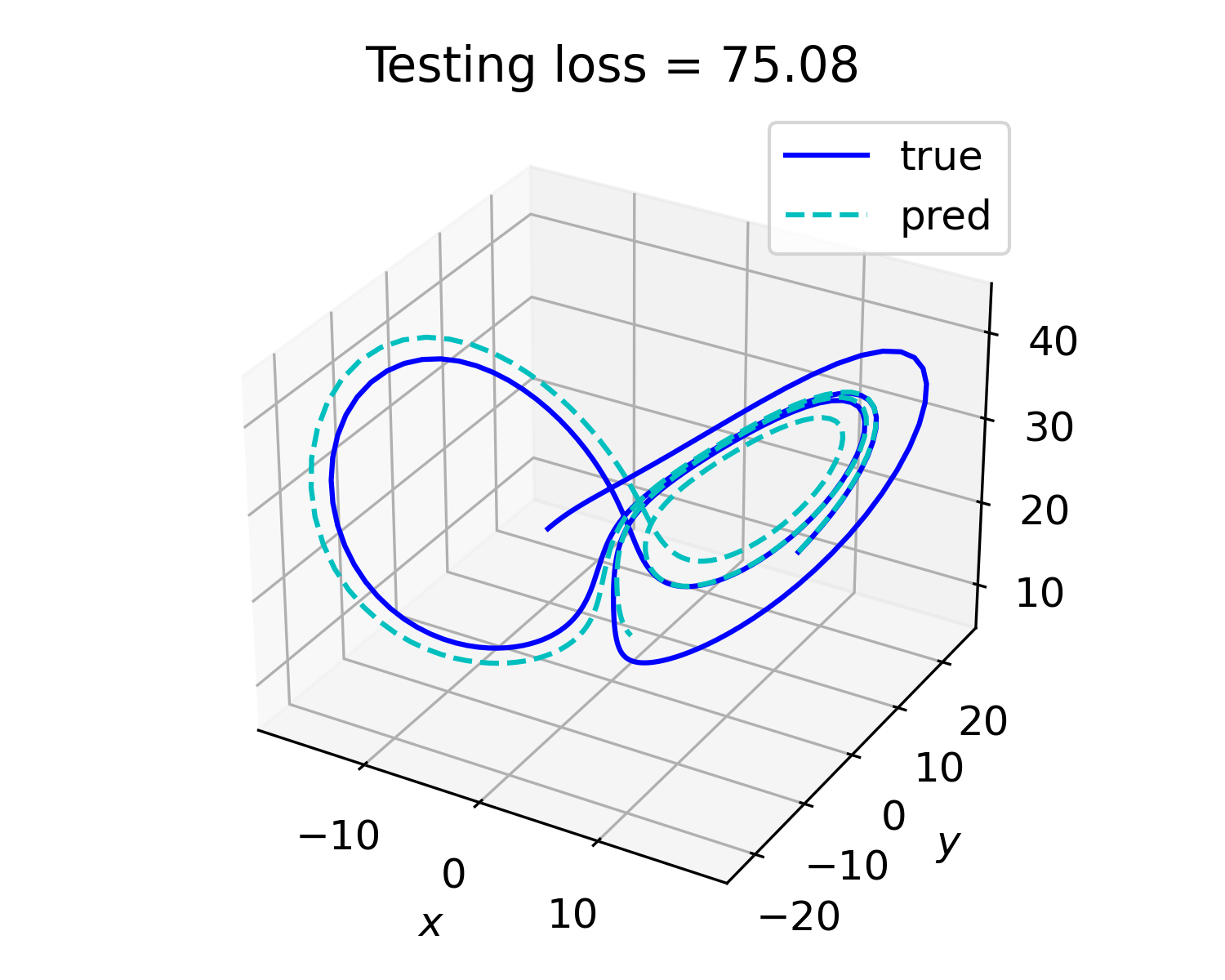}
\caption{Test (unseen)}
\end{subfigure}
\caption{Generalization performance of OCN on the Lorenz system. Solid lines represent the true trajectory, and the dashed lines represent the prediction given by OCN. }
\label{fig:ocn_pre}
\end{figure}

We first test the generalization performance of OCN by applying it to initial points that are different from the initial points used in training. Specifically, we consider a unit ball $S = \{u\;|\;\|u-x_0\|\leq 1\}$ where $x_0=[10, 15, 17]$, see Figure \ref{fig:ocn_pre} (a). The training data consists of 3 trajectories with the initial points in $S$, over time interval $[0, 3]$, and time step $\Delta t=0.01$. The training results are presented in Figure \ref{fig:ocn_pre} (b) (c) (d). We observe excellent agreements between the prediction by OCN and the true trajectories.

After training, we randomly select 300 points from $S$ as initial points. For each initial point, we generate the true trajectory data by (\ref{Lorenz}) and the prediction by OCN, then compute the loss using (\ref{loss}). The histogram of the testing loss over 300 trajectories is presented in Figure \ref{fig:ocn_pre} (e), from which we see that the testing loss is less than $80$ in over $80\%$ cases. In Figure \ref{fig:ocn_pre} (f) (g) (h), we present trajectories generated by $3$ different initial points, each corresponding to a different loss. Overall, OCN shows reasonably good prediction performance on data that is close to but does not belong to the training data.

\begin{figure}[ht]
\begin{subfigure}[b]{0.33\linewidth}
\centering
\includegraphics[width=1\linewidth]{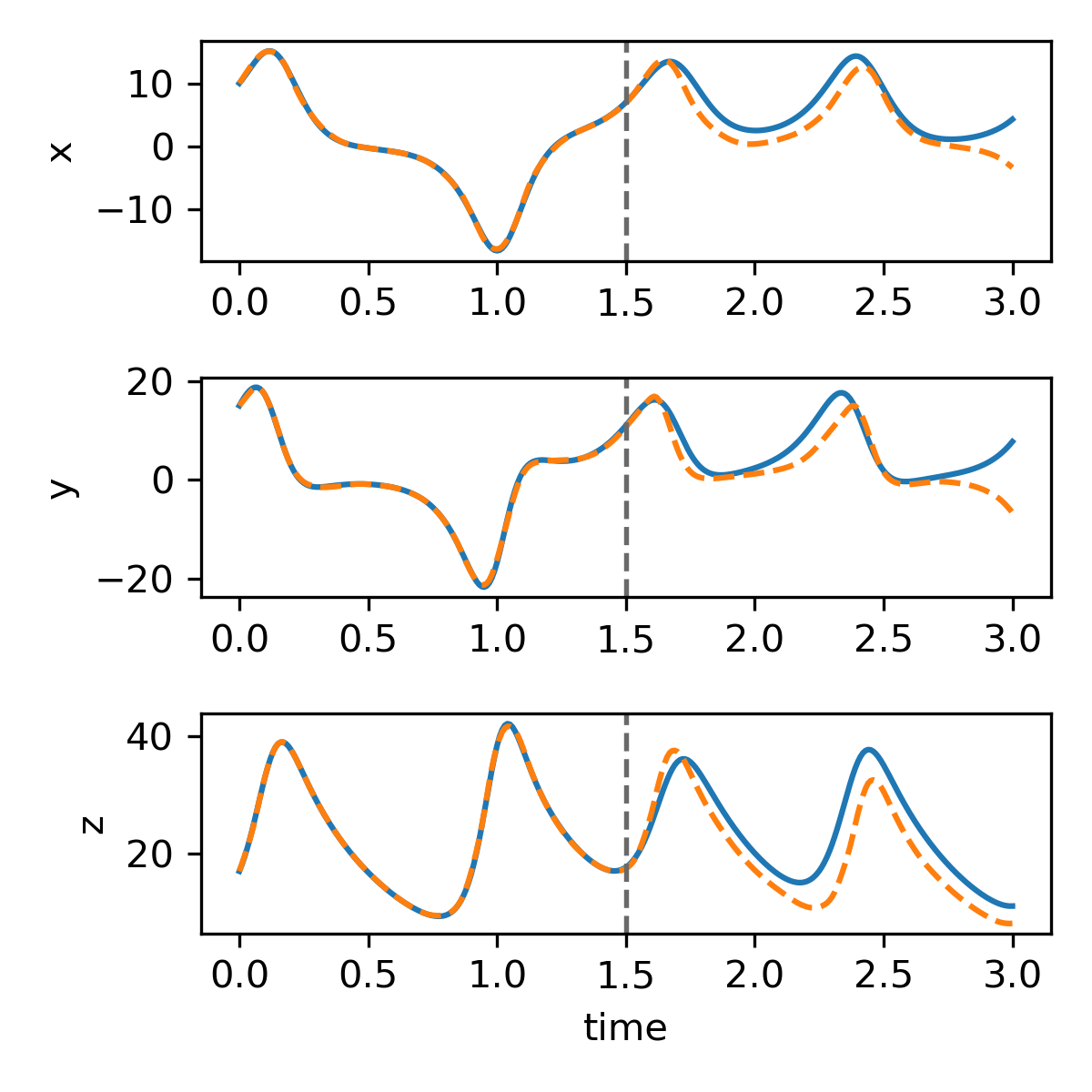}
\caption{OCN, $\Delta t=0.01$, no $\dot x$}
\end{subfigure}%
\begin{subfigure}[b]{0.33\linewidth}
\centering
\includegraphics[width=1\linewidth]{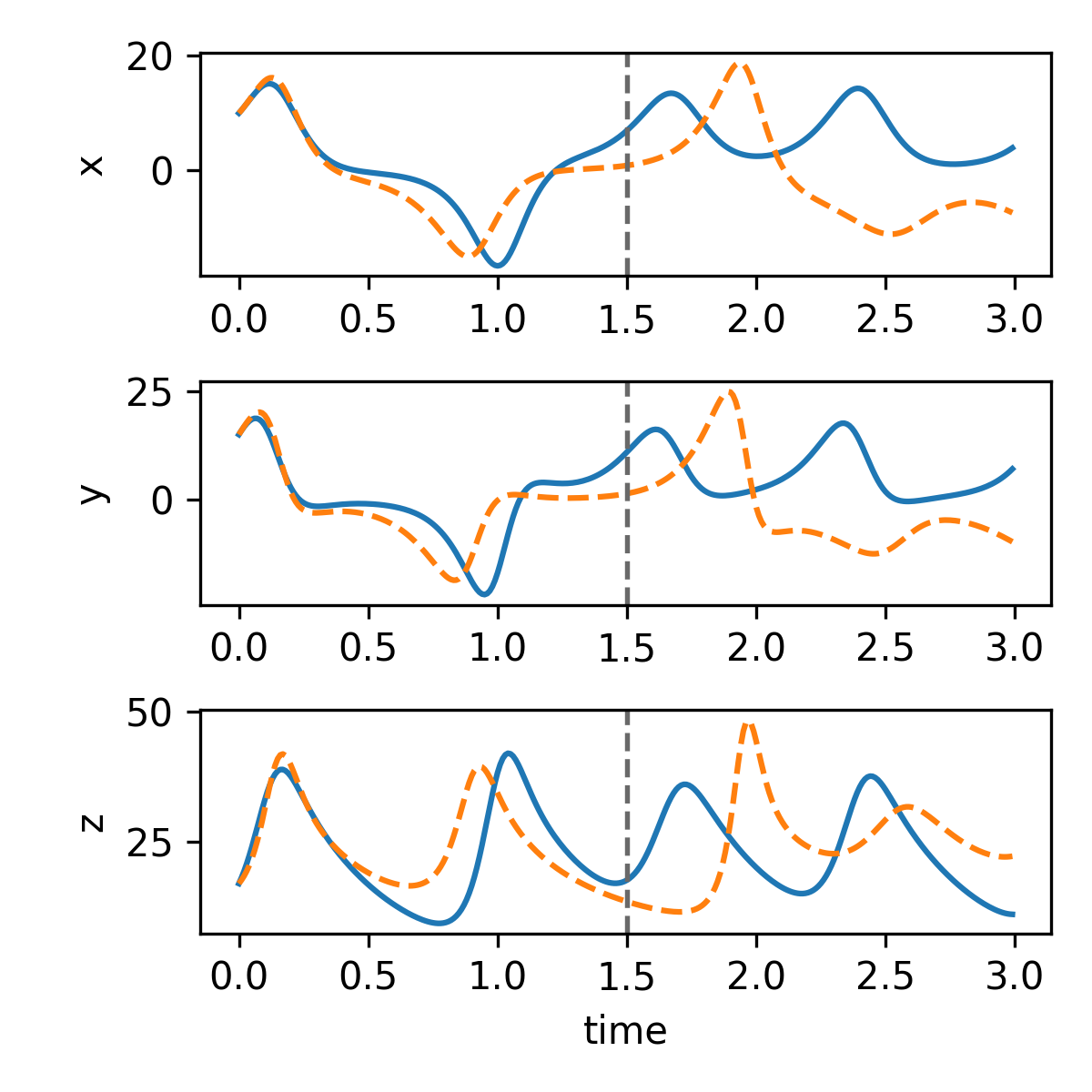}
\caption{SINDy, $\Delta t=0.001$, with $\dot x$}
\end{subfigure}%
\begin{subfigure}[b]{0.33\linewidth}
\centering
\includegraphics[width=1\linewidth]{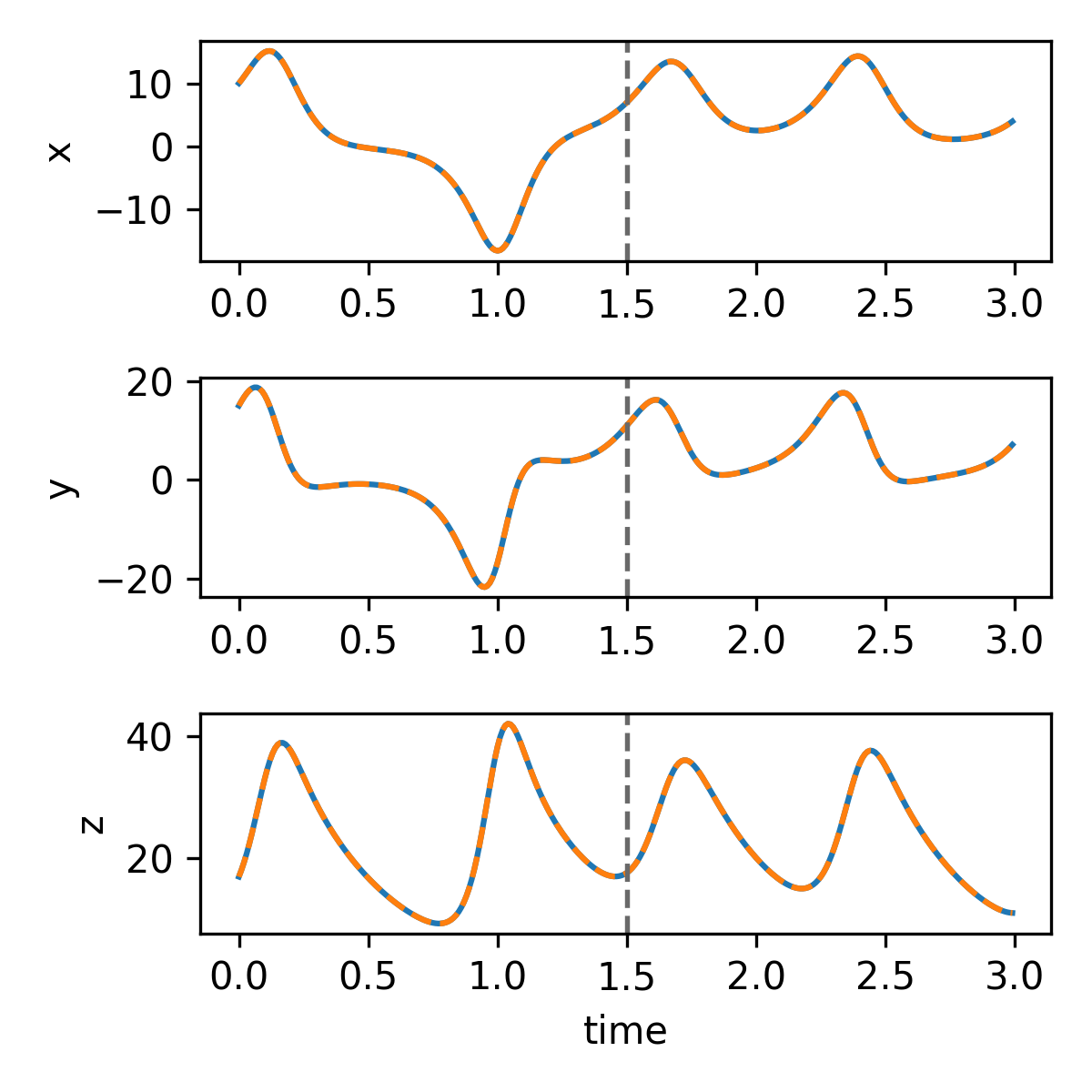}
\caption{SINDy, $\Delta t=0.00001$, no $\dot x$}
\end{subfigure}%
\caption{Comparison of the short-time performance between OCN and SIND on the Lorenz system. The results on $[0, 1.5]$ show the performance on the training data; the results on $[1.5, 3]$ show the prediction performance of each model.}
\label{fig:stc}
\end{figure}
For data-driven system discovery, the sparse identification of nonlinear dynamics (SINDy) \cite{BP+16} is a widely used method. It casts the system identification as a sparse regression problem over a large set of nonlinear library functions to find the fewest active terms that accurately reconstruct the system dynamics. The success of SINDy has inspired a large number of extensions and variants tailored for more specific problems  \cite{CL+19, RB+17, S17, ZA+18}. An obvious difference between SINDy and OCN is that SINDy, as its main feature, provides an explicit formula for the system, while OCN only gives network representations. Also, the derivative data $\dot x$ plays an important role in the framework of SINDy, while OCN does not require the information of $\dot x$.

In the next two subsections, we compare the performance of OCN with SINDy under two scenarios;  given short trajectory data or (relatively) long trajectory data. We consider different settings, including training data of different time steps $\Delta t$, with or without the derivative data $\dot x$.  When $\dot x$ is unavailable,  finite difference is used for SINDy to access estimations of $\dot x$. The comparison results are summarized in Table \ref{tb}.

\subsubsection{Short-time performance} In this case, the data used to train OCN is collected from $1$ trajectory with the initial point $[10, 15, 17]$, time interval $[0,1.5]$, and time step $\Delta t=0.01$. 
The training data for SINDy is collected from the same trajectory, while the time step is taken as $\Delta t = 0.001$. Also, the derivative data $\dot x$ is collected.

After the models are well trained, we apply them to generate trajectories on time interval $[0, 3]$ with the same initial point. The results are presented in Figure \ref{fig:stc}. We observe that compared with SINDy, OCN fits the data on $[0, 1.5]$ well and gives a good prediction on $[1.5, 3]$. The equation learned by SINDy is 
\begin{equation}
\begin{aligned}
&\dot x = 10(y-x),\\
&\dot y = x(28 - z) - y,\\
&\dot z = 0.034 z -0.091 z^2 + 0.034 x y z.
\end{aligned}    
\end{equation}
We see that (in this case), SINDy has difficulty in capturing the structure of the 3rd equation.

\subsubsection{Long-time performance} 
\begin{figure}[ht]
\begin{subfigure}[b]{0.5\linewidth}
\centering
\includegraphics[width=1\linewidth]{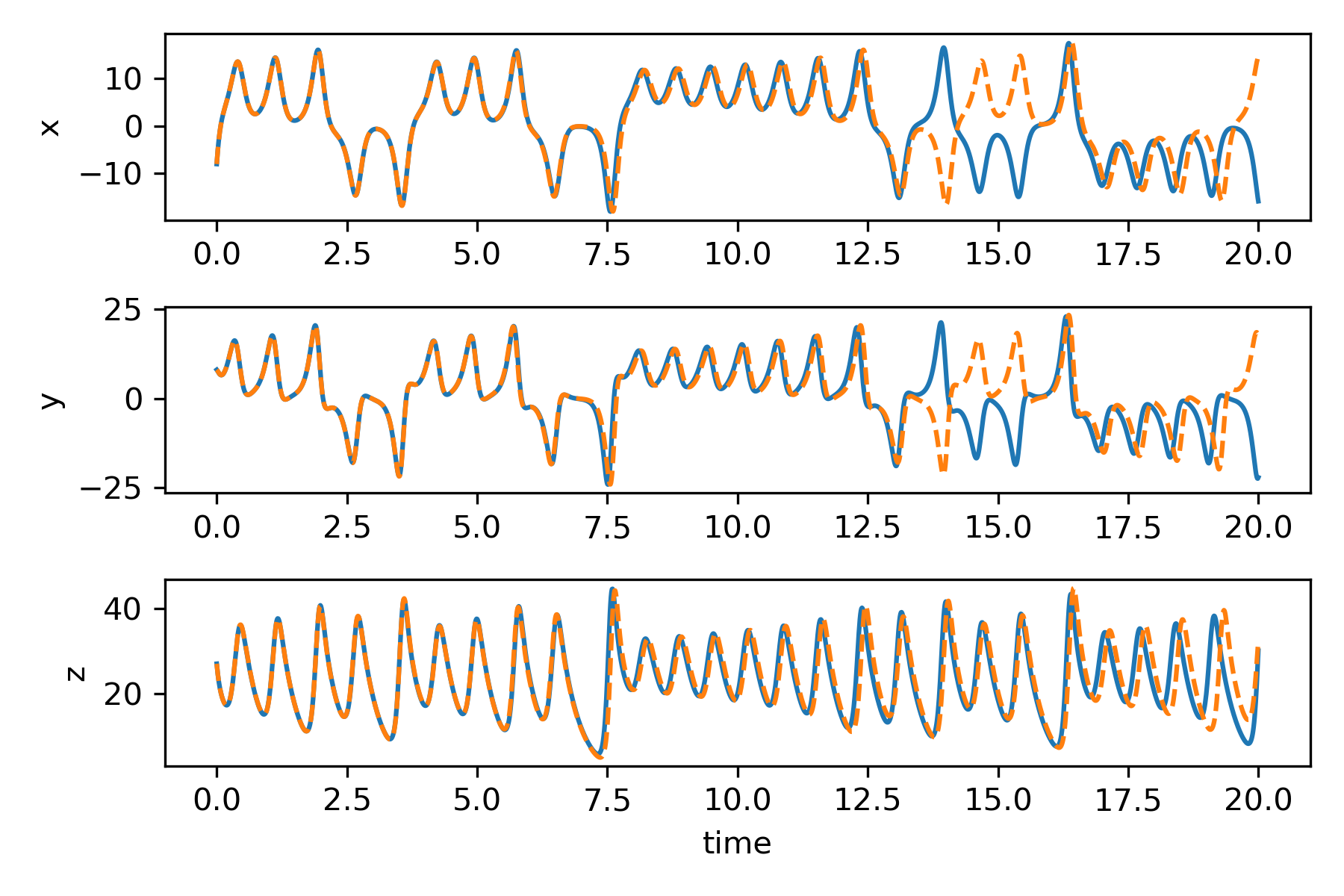}
\caption{OCN, $\Delta t=0.01$, no $\dot x$}
\end{subfigure}%
\begin{subfigure}[b]{0.5\linewidth}
\centering
\includegraphics[width=1\linewidth]{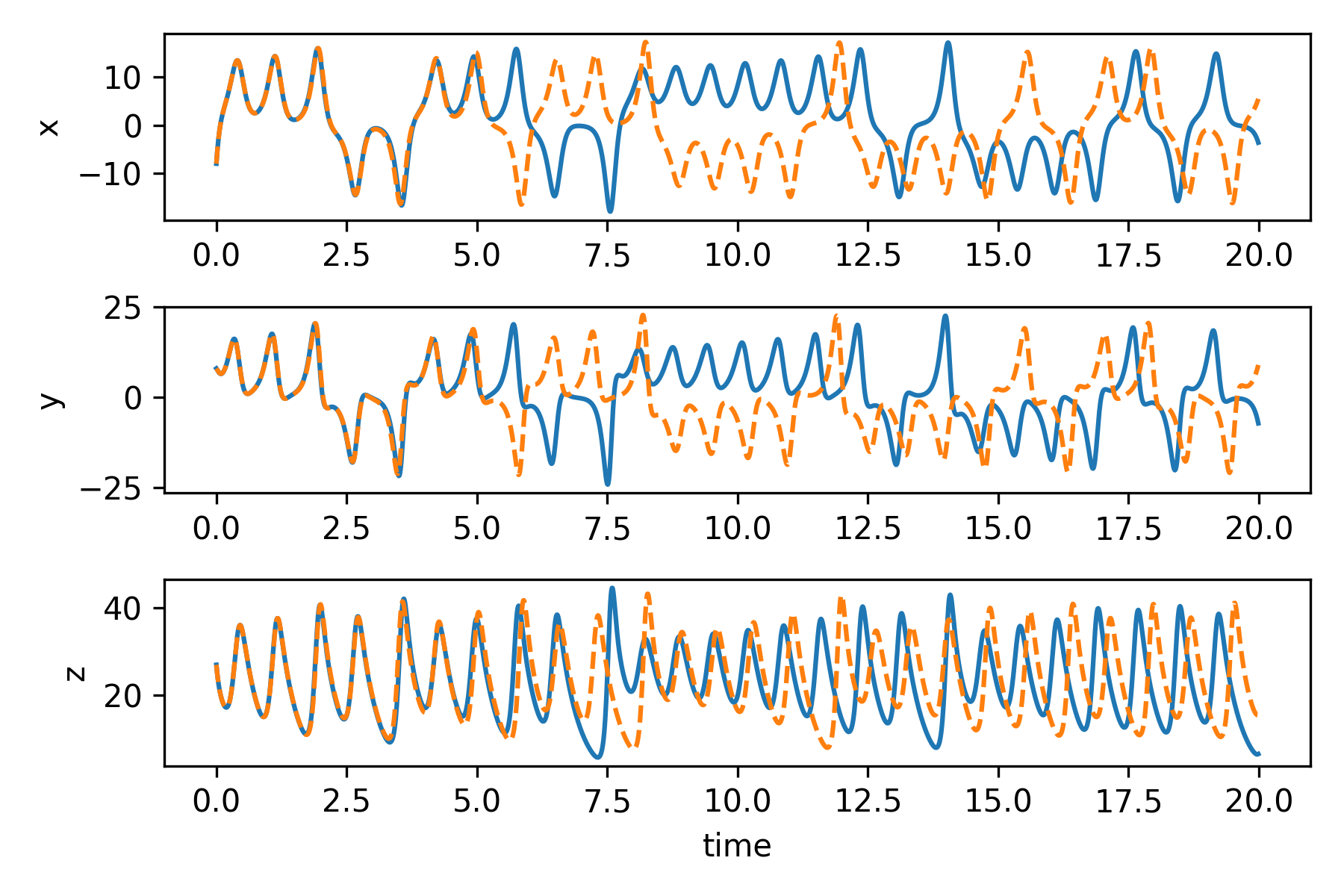}
\caption{SINDy, $\Delta t=0.01$, no $\dot x$}
\end{subfigure}%
\caption{Comparsion of the long-time performance between OCN and SINDy on the Lorenz system.}
\label{fig:l1}
\end{figure}

\begin{figure}[ht]
\begin{subfigure}[b]{0.5\linewidth}
\centering
\includegraphics[width=0.8\linewidth]{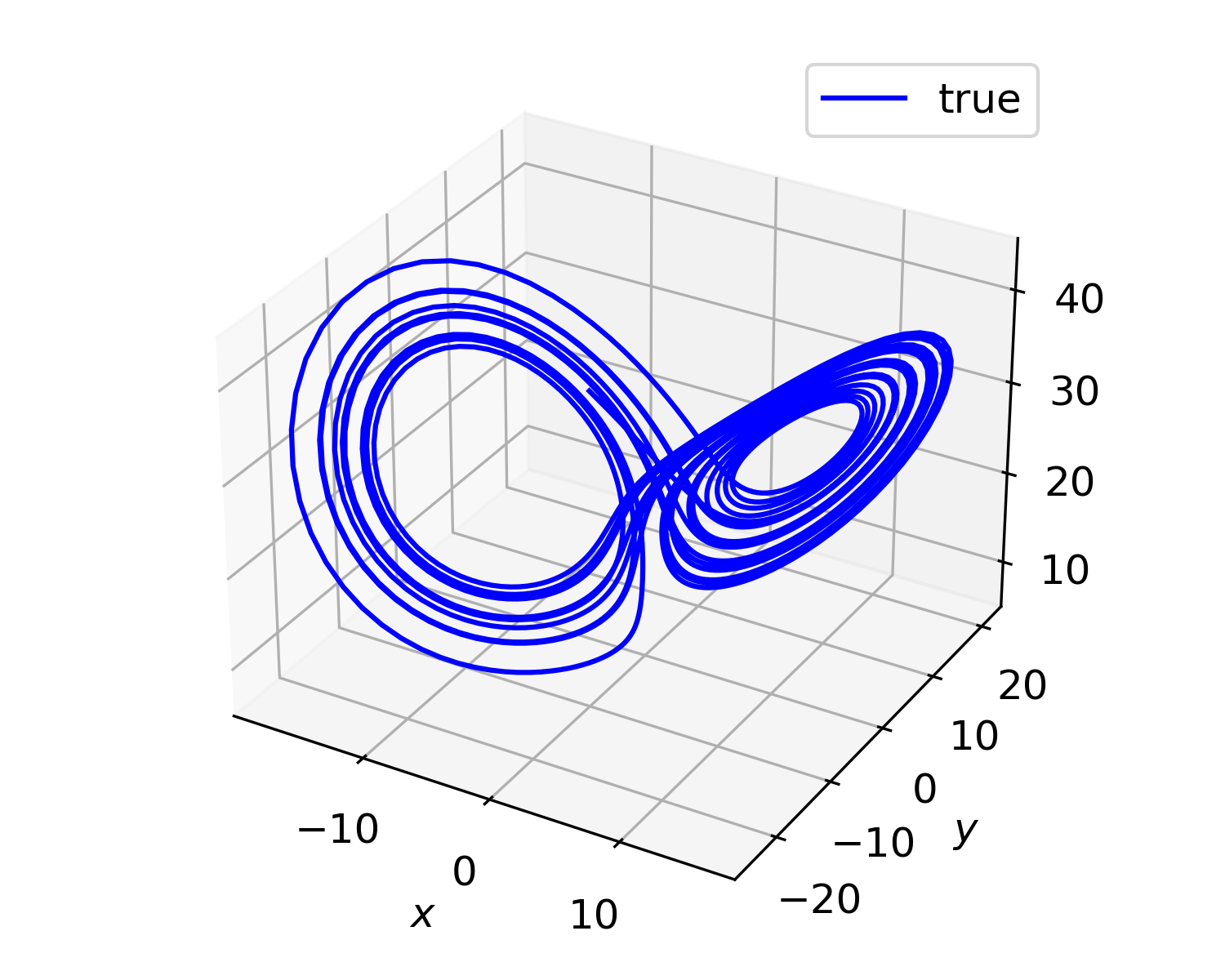}
\end{subfigure}%
\begin{subfigure}[b]{0.5\linewidth}
\centering
\includegraphics[width=0.8\linewidth]{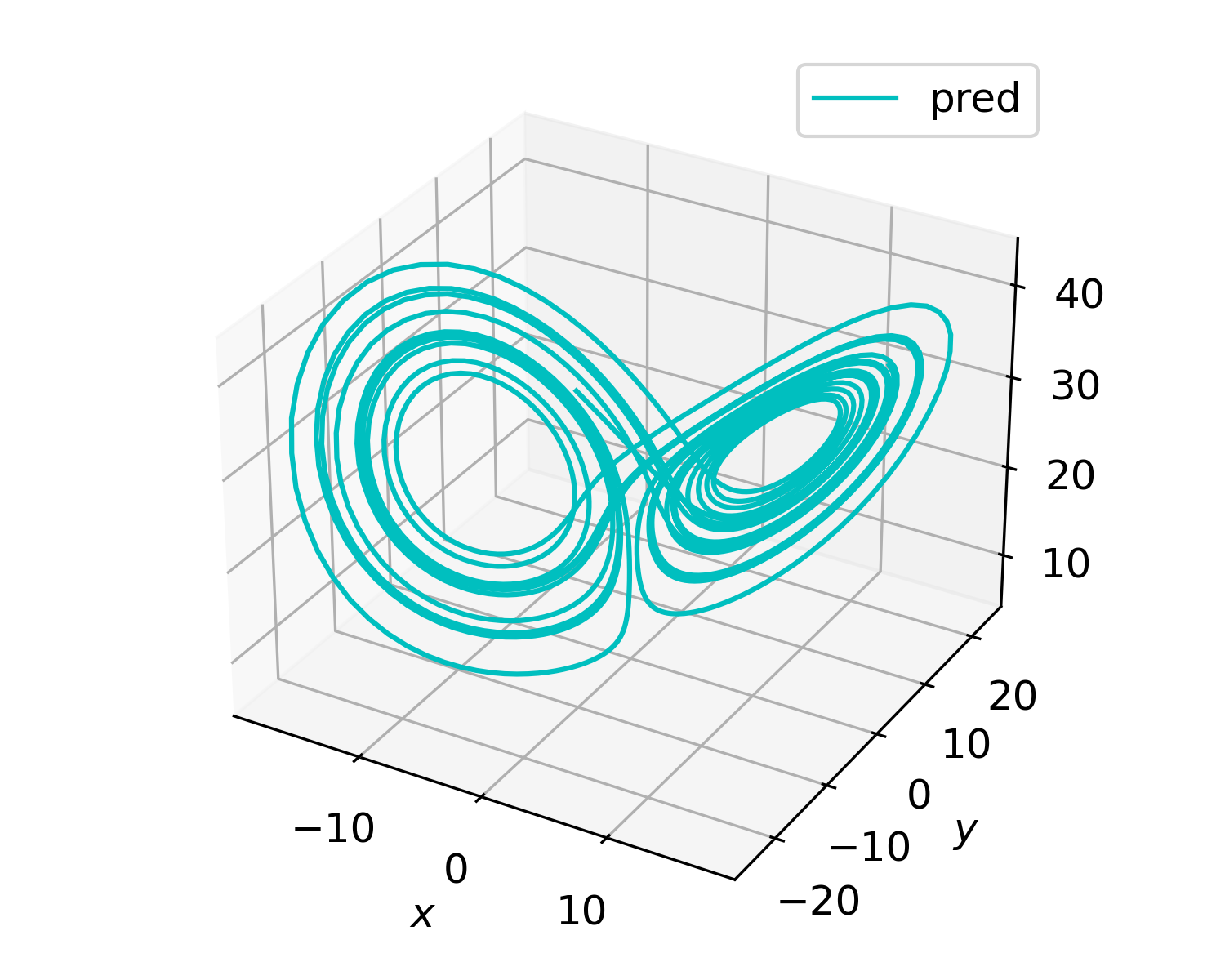}
\end{subfigure}%
\caption{True dynamics and prediction of OCN.}
\label{fig:l2}
\end{figure}

The data used for training OCN is collected from $1$ trajectory starting from $[-8, 8, 27]$, with time interval $[0,20]$ and time step $\Delta t=0.01$. The training results are presented in Figure \ref{fig:l1} (a) and Figure \ref{fig:l2}. The Lorenz system has a positive Lyapunov exponent, and small discrepancies between the true dynamics and learned models can grow exponentially, which should explain the large errors at a later time.  


These comparative assessments of neural network-based representation of dynamics versus an interpretable symbolic approach to representation suggest interesting tradeoffs between these choices for practitioners. 
Approaches like SINDy are simpler to implement, computationally more efficient in terms of model calibration, and interpretable. However, their performance relies very heavily on the accuracy of data $\dot x$. Moderately noisy $\dot x$ produces significant performance degradation. In contrast, OCNs are not interpretable, however, no data on $\dot x$ is required. Referring to the results in Table \ref{tb}, in cases $\Delta t$ is small e.g. $\Delta t=0.0001$, SINDy works very well. While in cases $\Delta t$ is relatively large e.g. $\Delta t= 0.01$, and without data on $\dot x$, OCN shows superior performance than SINDy, as also shown in Figure \ref{fig:l1}.
Overall, we find that in settings where (i) the observation data is collected from short-time trajectories, (ii) the derivative data $\dot x$ is unavailable, or (iii) the data $x$ has a relatively large time step $\Delta t$, OCN gives more accurate approximation than SINDy. A hybrid method that benefits from the advantages of the two approaches is certainly desirable; see e.g., \cite{SZS20,CL+19} for related works in this direction.

\begin{table} 
\begin{center}
\caption{Comparison of neural based (OCN) and symbolic regression (SINDy) approaches on the Lorenz system, using training data of different time steps $\Delta t$, with or without the derivative data $\dot x$. Training interval is the time interval from which the training data is collected. For cases with training interval $[0, 1.5]$, the loss is computed over time interval $[0, 3]$; for cases with training interval $[0, 20]$, the loss is computed over time interval $[0, 20]$.}
\label{tb} 
\begin{tabular}{c|c|p{5em}|p{5em}|c} 
\toprule
& Training interval & $\Delta t$ & $\dot x$ & Loss \\ 
\midrule
\multirow{1}{3em}{SINDy} 
& [0, 1.5] & 0.001 & yes & 108.23 \\ 
\midrule
OCN & [0, 1.5] & 0.01 & no & 6.93 \\ 
\midrule
\multirow{5}{3em}{SINDy} 
& [0, 20] & 0.01 & yes & 1.75e-6 \\ 
& [0, 20] & 0.01 & no & 124.96 \\ 
& [0, 20] & 0.001 & no & 55.48 \\ 
& [0, 20] & 0.0001 & no & 34.01 \\ 
\midrule
OCN & [0, 20] & 0.01 & no & 34.57 \\ 
\bottomrule
\end{tabular}
\end{center}
\end{table}

\section{Discussion}
This paper presents an approach to discovering gradient flows from data without assumptions on the form of the governing equations. We build on prior work in data-driven discovery of dynamical systems using machine learning techniques but with innovations related to a global network representation of the force field and an optimal control formulation, which allow our algorithm to scale to more complex problems.
The general form of the loss function allows for incorporating further knowledge of physics or regularization as necessary, so to make the method more accurate and robust. We derive error bounds for both the solution and the vector field. Specifically, we prove that the solution error depends on both the training error and the sparsity level of the time series data. We achieve this by carefully studying the error equation and obtaining a priori error bounds.

In numerical experiments, we demonstrate the effectiveness of OCN on a number of dynamical systems, including a linear gradient flow, a nonlinear gradient flow, the damped pendulum, and the chaotic Lorenz system. We show that OCN allows us to accurately learn the dynamics around different types of nodes, forecast future states, and maintain good generalization performance on testing data. Moreover, the comparison with SINDy on the chaotic Lorenz system illustrates the advantages of OCN when the data has a relatively large time step or the derivative data is not given. There are many dynamical systems to which this method may be applied, where there is ample data with the absence of governing equations. 

We see several avenues for future work, both theoretical and computational. For example, assuming the data is collected from the solution trajectory, we were able to improve the error bounds for $\|\nabla f - \nabla G\|$. What if we assume more structure on the dynamics? How can we improve the computational efficiency of solving the coupled control system? Can we deploy this to learn the dynamics of truly large-scale problems?

Let us also briefly discuss possible extensions of our method.  For systems with time dependence, such as 
$
\dot x = F(x,t),
$
for which we consider the augmented system
\begin{align*}
\dot x = F(x,u),\quad\dot u &= 1.
\end{align*}
For systems with physical parameters, 
$
\dot x = F(x,\mu),
$
then $\mu$ can be appended to the dynamics in the following way
\begin{align*}
\dot x = F(x,u),\quad \dot u &= 0.
\end{align*}
It is then possible to use neural networks to represent $F(x, u)$. Our results should be of broad interest to control and machine learning researchers using neural networks for learning and control.

Finally, we would like to point out that gradient flows in the form of partial differential equations (PDEs) can be reduced to ODE systems by the method of lines so that our method could be applied. In future work, we shall explore the learning of some important PDEs. 

\appendix
\renewcommand{\theequation}{\thesection.\arabic{equation}}

\section{Proof of Theorem \ref{thm1}}
The computation of the gradient of $J$ can be realized by the following recipe when $y=y(t;\theta)$ has been found to solve the following forward problem: 
\begin{equation}\label{forward}
\dot y(t)=-\partial_y G(y(t), \theta), \quad y(0)=x_0.
\end{equation}

(i) Build an augmented functional (associated Lagrangian) $\mathcal{L}$, a functional of independent variables  $\tilde y, p, \theta$ defined by
\begin{equation*}
\mathcal{L}(\tilde y,  p,\theta) = \sum_{i=1}^{n}L_i(\tilde y(t_i)) - \int_{0}^{T}(\dot{\tilde{y}}(t)+\nabla_{\tilde y} G(\tilde y(t),\theta))^\top  p(t)dt,
\end{equation*}
where $p$ is the Lagrange multiplier, and can be chosen freely. Taking $\tilde y= y$, 
we have
\begin{equation}\label{LJtheta}
\mathcal{L}(y, p,\theta) = \sum_{i=1}^{n}L_i( y(t_i)) = J (\theta).   
\end{equation}
In order to evaluate $\partial_\theta J$, we proceed to calculate the first variation of $\mathcal{L}(\tilde y, p,\theta)$ at $(y, \theta)$, defined by 
$$
\delta \mathcal{L}(y, p, \theta)
: =\lim_{\tau \to 0} \frac{
 \mathcal{L}(y+\tau \delta y, p, \theta+\tau \delta \theta)- \mathcal{L}(y, p, \theta)
}{\tau},
$$
from which we will see why $p$ should be chosen as in (\ref{00}).   

(ii) Defining the adjoint-state equations for $p$. By formal calculations, we obtain 
\begin{align*}
&\quad\delta \mathcal{L}(y, p, \theta)\\
& = \delta \sum_{i=1}^{n}\bigg(L_i(y(t_i))-\int_{t_{i-1}}^{t_i}\Big( \dot{y}(t)+\partial_y G( y(t),\theta)\Big)^\top p(t)dt\bigg)\\
& =  \sum_{i=1}^{n}\bigg(\delta y(t_i)^\top \partial_y L_i(y(t_i))-\int_{t_{i-1}}^{t_i}\Big(\delta \dot y(t)+\delta \partial_y G( y(t),\theta)\Big)^\top  p(t)dt\bigg)\\
& =\sum_{i=1}^{n}\bigg(\delta y(t_i)^\top \partial_y L_i(y(t_i))
-\delta y(t_i)^\top p(t^-_i)+ \delta y(t_{i-1})^\top p(t^+_{i-1})
\\
&\qquad  +\int_{t_{i-1}}^{t_i} (\delta y)^\top \dot p(t) - \Big(\nabla^2_y G( y(t),\theta) \delta y + \partial_\theta \partial_y G( y(t),\theta) \delta \theta \Big)^\top p(t)dt\bigg)\\ 
& = \delta y(T)^\top \Big(\partial_y L_n(y(T)) - p(T)\Big) + \delta y(0)^\top  p(0)
+ \sum_{i=1}^{n-1}\delta y(t_i)^\top \Big(\partial_y L_i(y(t_i)) - p(t^-_i) + p(t^+_i)\Big) \\
&\qquad + \sum_{i=1}^{n}\bigg(\int_{t_{i-1}}^{t_i} (\delta y)^\top \Big(\dot p(t) -\big(\nabla^2_y G( y(t),\theta)\big)^\top p(t)\Big) 
- (\delta \theta)^\top \Big(\big(\partial_\theta \partial_y G( y(t),\theta)\big)^\top p(t)\Big) dt\bigg),
\end{align*}
where we have used integration by parts, and regrouping of terms. Since $y(0)=x_0$ is fixed, $\delta y(0)=0$; if $p$ is taken to satisfy (\ref{00}), then
$$
\delta \mathcal{L}(y, p, \theta) = -(\delta \theta)^\top \int_0^T \big(\partial_\theta \partial_y G( y(t),\theta)\big)^\top p(t) dt.
$$
(iii) Computation of the gradient of $J$. Recall  (\ref{LJtheta}), the first variation of $J(\theta)$ is actually $\nabla J \cdot \delta \theta$, we thus conclude 
\begin{equation*}
\nabla J= -\int_0^T  \big(\partial_\theta \partial_y G(y(t),\theta)\big)^\top p(t)dt,    
\end{equation*}
as asserted in (\ref{grad-}). 

\section{Proof of Theorem \ref{thm2}}\label{pf-thm2}
It suffices to prove that the stated result holds for any $t\in[0,T]$. Without loss of generality, we assume 
$t\in I_i:=(t_i, t_{i+1}]$ for some $i\in\{0,1, ..., n-1\}$. Using the notation
$$
e_k(t) := \|y_k(t) - x(t)\|,
$$
where $y_k= y(t,\theta_k)$ and (\ref{ode-nn}), (\ref{gf}), we get
$$
\frac{d}{dt}e^2_k 
= 2 (y_k - x)\cdot \frac{d}{dt}(y_k - x) 
\leq 2 e_k \|\nabla f(x) - \partial_y G(y_k,\theta_k)\|,
$$
which is estimated by the Cauchy-Schwarz inequality. This further implies
\begin{equation}\label{dek}
\begin{aligned}
\dot e_k 
&\leq \|\nabla f(x) - \partial_y G(y_k,\theta_k)\| \\
&\leq \|\nabla f(x) - \nabla f(y_k)\| + \|\nabla f(y_k) - \partial_y G(y_k,\theta_k)\| \\
&\leq L_f e_k + R(y_k).    
\end{aligned}    
\end{equation}
Here we used the assumption that $\nabla f$ is $L_f$ Lipschitz continuous and the notation
$$
R(y_k(t)) := \|\nabla f(y_k(t)) - \partial_y G(y_k(t),\theta_k)\|.
$$
Rewriting (\ref{dek}) against an integrating factor $e^{-L_f t}$ we obtain 
$$
\frac{d}{dt}(e^{-L_f t} e_k(t)) \leq e^{-L_f t} R(y_k(t)).
$$
Integration of this over $(t_i,t)$ gives
\begin{equation}\label{ekbd}
\begin{aligned}
e_k(t) 
&\leq e^{L_f(t-t_i)} e_k(t_i) + \int_{t_i}^t e^{L_f(t-s)} R(y_k(s))ds\\
&\leq e^{L_f\Delta t} \Big(e_k(t_i) + \Delta t \max_{t\in I_i} R(y_k(t))\Big),
\end{aligned}    
\end{equation}
where $|t_{i+1}-t_i|\leq \max_i |t_{i+1}-t_i| =:\Delta t $ is used. 

We now proceed to bound the right hand side (RHS) of (\ref{ekbd}). First notice that
\begin{equation}\label{ekJ}
e_k(t_i) = \sqrt{\|y_k(t_i) - x_i\|^2} \leq \sqrt{J(\theta_k)}.
\end{equation}
For $R(y_k(t))$, we use triangle inequality to obtain
\begin{equation*}
\begin{aligned}
R(y_k(t)) 
&\leq \|\nabla f(y_k(t)) - \nabla f(y_k(t_i))\| \\
&\quad + \|\partial_y G(y_k(t_i),\theta_k) - \partial_y G(y_k(t),\theta_k)\| \\
&\quad + \|\nabla f(y_k(t_i)) - \partial_y G(y_k(t_i),\theta_k)\|, 
\end{aligned}
\end{equation*}
which implies
\begin{equation}\label{Ryk}
\max_{t\in I_i} R(y_k(t)) \leq D_1 + D_2 + D_3,    
\end{equation}
where
\begin{align*}
D_1 & = \max_{t\in I_i} \|\nabla f(y_k(t)) - \nabla f(y_k(t_i))\|, \\
D_2 & = \max_{t\in I_i} \|\partial_y G(y_k(t_i),\theta_k) - \partial_y G(y_k(t),\theta_k)\|, \\
D_3 & = \|\nabla f(y_k(t_i)) - \partial_y G(y_k(t_i),\theta_k)\|.
\end{align*}
We further derive bounds on $D_1, D_2, D_3$. The derivation of bounds on $D_1$ and $D_2$ are similar. The idea is to use $L_f$ Lipschitz continuity of $\nabla f$ and $L_{G_y}$, respectively with respect to $y$ to get
\begin{align*}
& D_1
\leq L_f \max_{t\in I_i} \|y_k(t) - y_k(t_i)\|,\\
& D_2
\leq L_{G_y} \max_{t\in I_i} \|y_k(t) - y_k(t_i)\|,
\end{align*}
then show the following bound
\begin{equation}\label{xybd}
\max_{t\in I_i} \|y_k(t) - y_k(t_i)\| 
\leq \frac{\Delta t}{1-\Delta t L_{G_y}} \Big(\|\partial_y G(x_i,\theta_k)\| + L_{G_y} \sqrt{J(\theta_k)}\Big).
\end{equation}
Hence for $\Delta t \leq \frac{1}{2L_{G_y}}$, we have
\begin{equation}\label{D12}
D_1 + D_2
\leq C_0\Delta t,
\end{equation}
where 
$$
C_0 = 2\Big(\|\partial_y G(x_i,\theta_k)\| + L_{G_y} \sqrt{J(\theta_k)}\Big) (L_f + L_{G_y}).
$$
For the derivation of (\ref{xybd}), we start with
\begin{equation}\label{maxdk}
\max_{t\in I_i} \|y_k(t) - y_k(t_i)\| 
= \max_{t\in I_i} \|\int_{t_i}^{t} \partial_y G(y_k(s),\theta_k)ds\| 
\leq \Delta t \max_{t\in I_i} \|\partial_y G(y_k(t),\theta_k)\|.    
\end{equation}
Using the $L_{G_y}$ Lipschitz continuity of $\partial_y G$ with respect to $y$, we have
\begin{equation*}
\begin{aligned}
\|\partial_y G(y_k(t),\theta_k) - \partial_y G(x_i,\theta_k)\| 
&\leq L_{G_y} \|y_k(t) - x_i\| \\
&\leq L_{G_y} (\|y_k(t) - y_k(t_i)\| + \|y_k(t_i) - x_i\|),
\end{aligned}    
\end{equation*}
which together with (\ref{ekJ}) lead to
\begin{equation}\label{maxG}
\max_{t\in I_i} \|\partial_y G(y_k(t),\theta_k)\| 
\leq \|\partial_y G(x_i,\theta_k)\| + L_{G_y} \sqrt{J(\theta_k)} + L_{G_y} \max_{t\in I_i} \|y_k(t) - y_k(t_i)\|.    
\end{equation}
Connecting (\ref{maxdk}) and (\ref{maxG}), we obtain (\ref{xybd}).

For the bound on $D_3$, we use triangle inequality to get
\begin{equation}\label{D3}
\begin{aligned}
D_3 
&\leq \|\nabla f(y_k(t_i)) - \nabla  f(x_i)\| + \|\nabla  f(x_i) + \frac{x_{i+1}-x_i}{\Delta t}\| \\
&\quad + \|-\frac{y_k(t_{i+1})-y_k(t_i)}{\Delta t} - \partial_y G(y_k(t_i),\theta_k)\|
+ \|\frac{y_k(t_{i+1})-y_k(t_i)}{\Delta t} -\frac{x_{i+1}-x_i}{\Delta t}\|.
\end{aligned}    
\end{equation}
The first term on the RHS of (\ref{D3}) can be bounded by
\begin{equation}\label{D31}
\|\nabla  f(y_k(t_i)) - \nabla f(x_i)\| 
\leq L_f e_k(t_i) \leq L_f \sqrt{J(\theta_k)},   
\end{equation}
using the $L_f$ Lipschitz continuous of $\nabla f$ and (\ref{ekJ}). 

For the second and third term on the RHS of (\ref{D3}), note that Assumption 1 and 2 also imply
\begin{align*}
& x(t_{i+1}) 
\leq x(t_i) - \Delta t \nabla f(x(t_i)) + \frac{L_f}{2} (\Delta t)^2,\\
& y(t_{i+1}) 
\leq y(t_i) - \Delta t \partial_y G(y_k(t_i),\theta_k) + \frac{L_{G_y}}{2} (\Delta t)^2.
\end{align*}
Since $x(t_i)=x_i$, we have
\begin{equation}\label{D324}
\begin{aligned}
& \|\nabla f(x_i) + \frac{x_{i+1}-x_i}{\Delta t}\| 
\leq \frac{L_f}{2} \Delta t,\\
& \|-\frac{y_k(t_{i+1})-y_k(t_i)}{\Delta t} - \partial_y G(y_k(t_i),\theta_k)\| 
\leq \frac{L_{G_y}}{2} \Delta t.
\end{aligned}    
\end{equation}

For the last term on the RHS of (\ref{D3}), we use triangle inequality and (\ref{ekJ}) to get
\begin{equation}\label{D33}
\begin{aligned}
&\quad\; \|\frac{y_k(t_{i+1})-y_k(t_i)}{\Delta t} - \frac{x_{i+1}-x_i}{\Delta t}\|\\
&\leq \frac{1}{\Delta t}\Big(\|y_k(t_{i+1}) - x_{i+1}\|+\|y_k(t_i) - x_i\|\Big) \leq \frac{2\sqrt{J(\theta_k)}}{\Delta t}.
\end{aligned}
\end{equation}
Substituting (\ref{D31}), (\ref{D324}), (\ref{D33}) into (\ref{D3}), we obtain the following bound on $D_3$
\begin{equation}\label{bd3}
D_3 \leq L_f \sqrt{J(\theta_k)} + \frac{L_f + L_{G_y}}{2} \Delta t + \frac{2\sqrt{J(\theta_k)}}{\Delta t}.
\end{equation}

With bounds on $D_1, D_2$ in (\ref{D12}) and $D_3$ in (\ref{bd3}), (\ref{Ryk}) becomes
$$
\max_{t\in I_i} R(y_k(t)) \leq C_0\Delta t + (L_f+\frac{2}{\Delta t})\sqrt{J(\theta_k)}.
$$
This together with (\ref{ekJ}), (\ref{ekbd}) and  $\Delta t \leq \frac{1}{2L_{G_y}}$ leads to
\begin{align*}
e_k(t) 
&\leq e^{L_f\Delta t} \Big(\sqrt{J(\theta_k)} + C_0(\Delta t)^2 + (L_f \Delta t+2)\sqrt{J(\theta_k)}\Big),\\
&\leq C_1\Big(\sqrt{J(\theta_k)} + (\Delta t)^2\Big),
\end{align*}
where 
$$
C_1 = e^{\frac{L_f}{2L_{G_y}}}\max\bigg\{C_0, 3+\frac{L_f}{2L_{G_y}}\bigg\},
$$
which further implies (\ref{41}) in Theorem \ref{thm2}.

The method used to derive (\ref{42}) is similar as that used for $D_3$. For any $i\in\{1,...,n\}$,
\begin{equation}\label{D3+}
\begin{aligned}
&\quad\;\|\nabla f(x_i) - \partial_y G(x_i,\theta_k)\|\\
&\leq \|\nabla  f(x_i) + \frac{x_{i+1}-x_i}{\Delta t}\| + \Big\|-\frac{y_k(t_{i+1})-y_k(t_i)}{\Delta t} - \partial_y G(y_k(t_i),\theta_k)\Big\| \\
&\quad + \|\nabla G(y_k(t_i)) - \nabla  G(x_i,\theta_k)\|+ \Big\|\frac{y_k(t_{i+1})-y_k(t_i)}{\Delta t} -\frac{x_{i+1}-x_i}{\Delta t}\Big\|.
\end{aligned}    
\end{equation}
Using (\ref{D324}), (\ref{D33}) and
\begin{equation}\label{gradG}
\|\nabla G(y_k(t_i)) - \nabla  G(x_i,\theta_k)\| 
\leq L_{G_y} e_k(t_i) \leq L_{G_y}\sqrt{J(\theta_k)},    
\end{equation}
we obtain
\begin{align*}
\|\nabla f(x_i) - \partial_y G(x_i,\theta_k)\| 
&\leq L_{G_y} \sqrt{J(\theta_k)} + \frac{L_f + L_{G_y}}{2} \Delta t + \frac{2\sqrt{J(\theta_k)}}{\Delta t},\\
&\leq \frac{5}{2}\frac{\sqrt{J(\theta_k)}}{\Delta t} + \frac{L_f + L_{G_y}}{2}\Delta t,\\
&\leq C_2 \bigg(\frac{\sqrt{J(\theta_k)}}{\Delta t} + \Delta t\bigg)
\end{align*}
where 
$$
C_2 = \max\bigg\{\frac{5}{2}, \frac{L_f+L_{G_y}}{2}\bigg\}.
$$
This further implies (\ref{42}) asserted in Theorem \ref{thm2}.

\section{Proof of Theorem \ref{thm2-v}}\label{pf-thm2-v}
The notations and techniques used in this proof are essentially the same as that used in the proof for Theorem \ref{thm2}. The only difference is the decomposition of the error on the gradient. More precisely, instead of (\ref{D3+}), now we have 
\begin{equation}\label{D3++}
\begin{aligned}
&\quad\;\|\nabla f(x_i) - \partial_y G(x_i,\theta_k)\|\\
&\leq \|\nabla  f(x_i) + \frac{x_{i+1}-x_i}{\Delta t}\| + \|-\frac{x_{i+1}-x_i}{\Delta t} - \partial_y G(y_k(t_i),\theta_k)\| \\
&\quad + \|\nabla G(y_k(t_i)) - \nabla  G(x_i,\theta_k)\|.
\end{aligned}    
\end{equation}
The second term on the right side is now part of the loss function, hence can be bounded by $\sqrt{J(\theta_k)}$. Recall (\ref{D324}) and (\ref{gradG}) for the bounds on the other two terms, we have
\begin{align*}
\|\nabla f(x_i) - \partial_y G(x_i,\theta_k)\| 
&\leq \frac{L_f}{2} \Delta t + \sqrt{J(\theta_k)} + L_{G_y} \sqrt{J(\theta_k)} ,\\
&\leq \frac{L_f}{2} \Delta t +  (L_{G_y}+1) \sqrt{J(\theta_k)} ,\\
&\leq C_2 (\sqrt{J(\theta_k)} + \Delta t)
\end{align*}
where 
$$
C_2 = \max\bigg\{\frac{L_f}{2}, L_{G_y}+1\bigg\}.
$$

\section{Proof of Theorem \ref{Sinv}}\label{pf-Sinv}
Taking gradient of $y$ in (\ref{arkf}) with respect to $y(0)$ gives
\begin{equation}\label{delta}
\begin{aligned}
\delta_{l+1} &= \delta_l + \tau_l \sum_{i=1}^{s}b_i d_{li}, \\
d_{li} &:= \frac{\partial g_{li}}{\partial y(0)} = -\partial^2_y G(y_{li},\theta)^\top \delta_{li}, \\
\delta_{li} &= \delta_l + \tau_l \sum_{j=1}^{s}a_{ij}d_{lj}.
\end{aligned}    
\end{equation}
That is, in the interval $[0, T]$, $\delta(t)$ is discretized by the same method as $y(t)$. In each time interval $(t_{i-1}, t_i]$, we have 
\begin{align}\notag
\Delta 
&= \delta_{l+1}^\top p_{l+1} - \delta_{l}^\top p_{l} \\\notag
&= \bigg(\delta_l + \tau_l \sum_{i=1}^{s}b_i d_{li}\bigg)^\top \bigg(p_{l} + \tau_l \sum_{i=1}^{s}\tilde b_i h_{li}\bigg) - \delta_{l}^\top p_{l} \\\label{I123}
&= I_1 + I_2 + I_3,
\end{align}
where 
\begin{align*} 
I_1 = \tau_l \sum_{i=1}^{s}b_i d_{li}^\top p_l, \quad
I_2 = \tau_l \sum_{i=1}^{s}\tilde b_i \delta_l^\top h_{li}, \quad
I_3 = \tau_l^2 \sum_{i=1}^{s}\sum_{j=1}^{s} b_i\tilde b_j d_{li}^\top h_{li}.
\end{align*}
Below we deal with $I_1, I_2, I_3$ separately. 

For $I_1$, we note that if $b_i\neq 0$, then
\begin{align*}
p_{li} 
& = p_{l+1} - \tau_l \sum_{j=1}^{s}\tilde b_j\frac{a_{ji}}{b_i}h_{lj}  \\
& = p_{l} + \tau_l \sum_{i=1}^{s}\tilde b_i h_{li}  - \tau_l \sum_{j=1}^{s}\tilde b_j\frac{a_{ji}}{b_i}h_{lj} \\
& = p_{l} +  \tau_l \sum_{j=1}^{s}\tilde b_j\Big(1-\frac{a_{ji}}{b_i}\Big)h_{lj}.
\end{align*}
Denote $Q = \{i\;|\;b_i=0\}$, then $I_1$ can be rewritten as
\begin{align*}
I_1 
&= \tau_l \sum_{i\notin Q} b_i d_{li}^\top p_l \\
&= \tau_l \sum_{i\notin Q} b_i d_{li}^\top \bigg(p_{li} - \tau_l \sum_{j=1}^{s}\tilde b_j\Big(1-\frac{a_{ji}}{b_i}\Big)h_{lj}\bigg) \\
&= \tau_l \sum_{i\notin Q} b_i d_{li}^\top p_{li} 
- \tau_l^2 \sum_{i\notin Q}\sum_{j=1}^{s}(b_i\tilde b_j- \tilde b_j a_{ji})d_{li}^\top h_{lj}. 
\end{align*}
For $I_2$, we have
\begin{align*}
I_2
&= \tau_l \sum_{i=1}^{s}\tilde b_i \delta_l^\top h_{li} \\
&= \tau_l \sum_{i=1}^{s}\tilde b_i \Big(\delta_{li} - \tau_l \sum_{j=1}^{s}a_{ij}d_{lj}\Big)^\top h_{li} \\
&= \tau_l \sum_{i=1}^{s}\tilde b_i \delta_{li}^\top h_{li}
- \tau_l^2 \sum_{i=1}^{s}\sum_{j=1}^{s}\tilde b_i a_{ij}d_{lj}^\top h_{li}\\
&= \tau_l \sum_{i=1}^{s}\tilde b_i \delta_{li}^\top h_{li}
- \tau_l^2 \sum_{j=1}^{s}\sum_{i=1}^{s}\tilde b_j a_{ji}d_{li}^\top h_{lj}\\
&= \tau_l \sum_{i\notin Q} b_i \delta_{li}^\top h_{li} + \tau_l \sum_{i\in Q}\tilde b_i \delta_{li}^\top h_{li} \\
&\quad - \tau_l^2 \sum_{i\notin Q}\sum_{j=1}^{s} b_j a_{ji}d_{li}^\top h_{lj}
- \tau_l^2 \sum_{i\in Q}\sum_{j=1}^{s}\tilde b_j a_{ji}d_{li}^\top h_{lj}.
\end{align*}
Here $\tilde b_i = b_i$ for $i\notin Q$ was used in the last equality. 

For $I_3$, using the notation of $Q$, we have
\begin{equation*}
I_3 = \tau_l^2 \sum_{i\notin Q}\sum_{j=1}^{s} b_i\tilde b_j d_{li}^\top h_{li}.   
\end{equation*}

Addind up $I_1, I_2, I_3$, we can simplify (\ref{I123}) as 
\begin{align*}
\Delta =  \tau_l \sum_{i\notin Q} b_i (d_{li}^\top p_{li} + \delta_{li}^\top h_{li}) + \tau_l \sum_{i\in Q}\tilde b_i \delta_{li}^\top h_{li} - \tau_l^2 \sum_{i\in Q}\sum_{j=1}^{s}\tilde b_j a_{ji}d_{li}^\top h_{lj}.
\end{align*}
Note that for $i\in Q$, we have $\tilde b_i=\tau_l$ and $p_{li}  = - \sum_{j=1}^{s}\tilde b_ja_{ji}h_{lj}$, thus $\Delta$ can be further reduced as
\begin{equation*}
\Delta =  \tau_l \sum_{i\notin Q} b_i (d_{li}^\top p_{li} + \delta_{li}^\top h_{li})
+ \tau_l^2 \sum_{i\in Q} (d_{li}^\top p_{li} + \delta_{li}^\top h_{li}).
\end{equation*}
Now it suffices to show that 
\begin{equation}\label{disdh}
d_{li}^\top p_{li} + \delta_{li}^\top h_{li} = 0,\quad \forall l, i.
\end{equation}
This can be derived from the property that $\delta(t)^\top p(t)$ is conserved in the continuous level. Using (\ref{Sp}) and the notations $d_{li} = -\partial^2_y G(y_{li},\theta)^\top \delta_{li}$, $h_{li} = \partial^2_y G(y_{li},\theta)^\top p_{li}$, (\ref{disdh}) follows.

\bibliography{ref}
\bibliographystyle{amsplain}

\end{document}